\documentclass[11pt]{article}

\usepackage[all]{xy}
\usepackage{amsmath}
\usepackage{amsthm}
\usepackage{amscd}
\usepackage{amssymb}
\usepackage{amsopn}

\input{epsf.tex}

\newcounter{enumitemp}
\newenvironment{enumeratecontinue}{
  \setcounter{enumitemp}{\value{enumi}}
  \begin{enumerate}
  \setcounter{enumi}{\value{enumitemp}}
}
{
  \end{enumerate}
}

\numberwithin{equation}{section}

\newcommand{\nb}[1]{#1\nobreakdash-}

\newcommand\ds\displaystyle

\theoremstyle{definition}

\theoremstyle{plain}
\newtheorem*{theorem*}{Theorem}
\newtheorem{theorem}{Theorem}
\newtheorem{proposition}[theorem]{Proposition}
\newtheorem{lemma}[theorem]{Lemma}
\newtheorem{corollary}[theorem]{Corollary}

\newtheorem{fact}[theorem]{Fact}

\DeclareMathOperator{\Out}{Out}
\DeclareMathOperator{\Inn}{Inn}
\DeclareMathOperator{\Aut}{Aut}

\DeclareMathOperator{\Length}{Length}

\DeclareMathOperator{\Isom}{Isom}
\DeclareMathOperator{\Stab}{Stab}

\DeclareMathOperator\Edges{\mathcal{E}}

\newcommand\reals{{\mathbf R}}
\newcommand\R\reals
\newcommand\real\reals

\newcommand\Z{{\mathbf Z}}

\newcommand\inject{\hookrightarrow}

\newcommand\infinity{\infty}
\newcommand{\bdy}{\partial}
\newcommand{\from}{\colon}
\newcommand\composed{\circ}
\newcommand\suchthat{\bigm|}
\newcommand\inv{{-1}}
\newcommand\union{\cup}
\newcommand\abs[1]{\left| #1 \right|}
\newcommand\norm[1]{\left\| #1 \right\|}

\newcommand\intersect{\cap}
\newcommand\restrict{\bigm|}

\newcommand\Teichmuller{Teichm\"uller}

\newcommand\MCG{{\mathcal M \mathcal C \mathcal G}}

\newcommand\cross{\times}

\newcommand\C{\mathcal C}

\newcommand\E{\mathcal E}

\renewcommand\L{\mathcal L}

\renewcommand\P{\mathcal P}

\newcommand\X{\mathcal X}
\newcommand\<\langle
\renewcommand\>\rangle
\newcommand\wt\widetilde
\newcommand\wh\widehat
\newcommand\disjunion\sqcup
\newcommand\fol{{\cal F}}

\DeclareMathOperator\interior{int}

\newcommand\act\circlearrowright

\newcommand\Hull{{\mathcal H}}

\newcommand\BH{\cite{BestvinaHandel:tt}}
\newcommand\BookZero{\cite{BFH:laminations}}
\newcommand\BookOne{\cite{BFH:TitsOne}}
\newcommand\BookTwo{\cite{BFH:TitsTwo}}
\newcommand\OuterLimits{\cite{BestvinaFeighn:OuterLimits}}
\newcommand\StableActions{\cite{BestvinaFeighn:stableactions}}

\DeclareMathOperator\PF{PF}

\newcommand\geodesic[2]{\overline{#1 \, #2}}
\newcommand\bx{\textbf{x}}

\newcommand\TG{{\mathcal{TG}}}

\title{Parageometric outer automorphisms of free groups}
\author{Michael Handel\thanks{Supported in part by NSF grant DMS0103435.}\ \ 
and Lee Mosher\thanks{Supported in part by NSF grant DMS0103208.}}

\begin{document}
\maketitle

\begin{abstract} We study those fully irreducible outer automorphisms $\phi$ of a finite rank free group $F_r$ which are \emph{parageometric}, meaning that the attracting fixed point of~$\phi$ in the boundary of outer space is a geometric $\reals$-tree with respect to the action of~$F_r$, but~$\phi$ itself is not a geometric outer automorphism in that it is not represented by a homeomorphism of a surface. Our main result shows that the expansion factor of $\phi$ is strictly larger than the expansion
factor of $\phi^\inv$. As corollaries (proved independently by Guirardel), the inverse of a parageometric outer automorphism is neither geometric nor parageometric, and a fully irreducible outer automorphism $\phi$ is geometric if and only if its attracting and repelling fixed points in the boundary of outer space are geometric $\reals$-trees.
\end{abstract}

\section{Introduction}

There is a growing dictionary of analogies between theorems about the mapping class group of a surface $\MCG(S)$ and theorems about the outer automorphism group of a free group $\Out(F_r)$. For example, the Tits alternative for $\MCG(S)$ \cite{McCarthy:Tits} is proved using Thurston's theory of measured geodesic laminations \cite{FLP}, and for $\Out(F_r)$ it is proved using the Bestvina--Feighn--Handel theory of
laminations (\BookZero, \BookOne, \BookTwo).

\paragraph{Expansion factors.} Here is a result about $\MCG(S)$ of which one might hope to have an analogue
in $\Out(F_r)$. Given a finitely generated group $G$, its outer automorphism group $\Out(G)$ acts on the
set of conjugacy classes $\C$ of $G$. Given $c\in\C$ let $\norm{c}$ be the smallest word length of a
representative of $c$. Given $\phi \in \Out(G)$ define the \emph{expansion factor}
$$\lambda(\phi) = \sup_{c \in \C} \left( \limsup_{n \to +\infinity} \norm{\phi^n(c)}^{1/n}
\right)
$$
If $\phi \in \MCG(S) \approx \Out(\pi_1 S)$ is pseudo-Anosov then
$\lambda(\phi)$ equals the pseudo-Anosov expansion factor \cite{FLP}. By
combining results of Thurston and Bers one obtains:

\begin{theorem*} If $\phi\in\MCG(S)$ is pseudo-Anosov, then there is a unique $\phi$-invariant geodesic
in \Teichmuller\ space, consisting of the points in \Teichmuller\ space which minimize the translation
distance under $\phi$. This translation distance equals $\log(\lambda(\phi))$.
\end{theorem*}

Consider a free group $F_r$ and the Culler-Vogtmann \emph{outer space} $\X_r$ on which $\Out(F_r)$ acts
properly. Recall that $\phi \in \Out(F_r)$ is \emph{reducible} if there is a nontrivial free
decomposition $F_r = A_1 * \cdots * A_k * B$ such that $\phi$ permutes the conjugacy classes of
$A_1,\ldots,A_k$; otherwise, $\phi$ is irreducible. If $\phi^k$ is irreducible for all $k\ge 1$ then we
say that $\phi$ is \emph{fully irreducible}. By analogy, a mapping class on $S$ is reducible if it
preserves the isotopy classes of some nontrivial decomposition of $S$ into essential subsurfaces, and is
fully irreducible if and only if it is represented by a pseudo-Anosov homeomorphism.

Question: Does the above theorem have an analogue for a fully irreducible
$\phi \in \Out(F_r)$? The question does not quite make sense because no metric is specified on
$\X_r$, but one can instead ask: Is there an analogue with respect to \emph{some}
$\Out(F_r)$-equivariant metric on $\X_r$? Or on any other metric space on which $\Out(F_r)$ acts?

Answer: No. If translation distance for $\phi$ is uniquely minimized on an axis $\gamma$ then, by
symmetry of the distance function, translation distance for $\phi^\inv$ is also uniquely minimized on
$\gamma$, and the minima for $\phi$ and for $\phi^\inv$ along $\gamma$ are equal. By equating minimal
translation distance with $\log(\lambda)$, one would conclude that $\lambda(\phi) =\lambda(\phi^\inv)$.
However, an example from \BH\ has the property that
$\lambda(\phi)\ne\lambda(\phi^\inv)$: consider $\phi \in\Out(F_3)$ and $\phi^\inv$ represented by the
automorphisms
$$
\Phi \colon \begin{cases}
A &\to AC \\
B &\to A \\
C &\to B
\end{cases}
\qquad\qquad
\Phi^\inv \colon \begin{cases}
A &\to B \\
B &\to C \\
C &\to \overline B A
\end{cases}
$$
Interpreting these formulas as self maps of the three-petaled rose, each is clearly a
train track map. From the results of \BH\ (see also Proposition~\ref{PropStretchFactor}), 
if $f \from G \to G$ is a train track representative of $\phi \in \Out(F_r)$, and if $\PF(M_f)$
denotes the Perron-Frobenius eigenvalue of the transition matrix $M_f$,
then $\lambda(\phi) = \PF(M_f)$. For the above two train track maps we therefore obtain
$$\lambda(\phi) = \PF\begin{pmatrix}
1 & 0 & 1 \\
1 & 0 & 0 \\
0 & 1 & 0 \\
\end{pmatrix} > 1.4, 
\qquad\qquad \lambda(\phi^\inv) = \PF\begin{pmatrix}
0 & 1 & 0 \\
0 & 0 & 1 \\
1 & 1 & 0
\end{pmatrix} < 1.4
$$

Confronted with such a strange phenomenon, one strategy is to see what appropriately weaker results can
be proved. We follow this strategy in the companion paper \cite{HandelMosher:logratio} where we show
that the ratio $\log(\lambda(\phi)) / \log(\lambda(\phi^\inv))$ is bounded by a constant depending only
on the rank $r$. This is what one would expect if there were an axis $\gamma$ for $\phi$ with
translation distance $\log(\lambda(\phi))$ and an axis $\gamma'$ for $\phi^\inv$ with translation
distance $\log(\lambda(\phi^\inv))$, such that $\gamma,\gamma'$ are fellow travelers. Encouraged
by this result, we have pursued the study of axes in outer space, with some interesting analogues
of uniqueness of axes \cite{HandelMosher:axes}.  

\paragraph{Parageometric outer automorphisms.} In this work we pursue another strategy: explore the
strange phenomenon on its own terms. The $\phi$ described above turns out to be an example of a
\emph{parageometric} outer automorphism, as we discovered by comparing discussions of this same example
in \StableActions\ and in \OuterLimits. The interest in this concept was pointed out in
\cite{GJLL:index}, where we found the terminology ``parageometric''. 

While we believe that the phenomenon $\lambda(\phi) \ne \lambda(\phi^\inv)$ is generic among fully irreducible outer automorphisms $\phi$, we shall show that inequality \emph{always} holds when $\phi$ is parageometric, in fact we give an explicit geometric argument which shows that $\lambda(\phi) > \lambda(\phi^\inv)$. 

To define parageometricity, recall the action of $\Out(F_r)$ on the compactified outer space $\overline
\X_r = \X_r \union \bdy\X_r$ consisting of (classes of) very small actions of $F_r$ on $\reals$-trees
(see \cite{CohenLustig:verysmall} for simplicial $\reals$-trees and \OuterLimits\ for nonsimplicial $\reals$-trees); we shall call these objects ``$F_r$-trees''. The action of a fully irreducible $\phi\in \Out(F_r)$ on $\overline \X_r$ has source--sink dynamics, with a repelling $F_r$-tree $T_-\in\bdy \X_r$, and an attracting $F_r$-tree $T_+ \in \bdy\X_r$ (see \BookZero\ for orbits in $\X_r$ and \cite{LevittLustig:NorthSouth} for orbits in $\bdy\X_r$). An $F_r$-tree is \emph{geometric} if it is dual in the appropriate sense to a measured foliation defined on some \nb{2}complex whose fundamental group surjects to $F_r$ \cite{LevittPaulin:geometric}. For example, if a fully irreducible $\phi\in\Out(F_r)$ is geometric, meaning that it is represented by an automorphism of a surface with boundary, then both of the $F_r$-trees $T_-,T_+$ are geometric: this follows from Thurston's theorem that $\phi$ is represented by a pseudo-Anosov surface homeomorphism $f\from S\to S$, because $T_-$ and $T_+$ are dual to the stable and unstable measured foliations of $f$ defined on the surface $S$.

A fully irreducible $\phi \in \Out(F_r)$ is said to be \emph{parageometric} if $T_+$ is a geometric
$F_r$-tree but $\phi$ is not a geometric outer automorphism. 

For example, the outer automorphism $\phi$ described above is parageometric: geometricity of the
$F_r$-tree $T_+$ is proved in Example 3.4 of \OuterLimits; and Levitt's ``thinness'' property for
$T_-$ is proved in Example 10.1 of \StableActions, showing that $T_-$ is not a geometric $F_r$-tree,
and so $\phi$ is not a geometric outer automorphism. In Proposition~\ref{PropGeometric} we will gather
results of \OuterLimits\ and \BH\ which give a method of characterizing parageometricity solely from the
properties of a train track representative.

Here are our main results:

\begin{theorem}
\label{TheoremLambdaBigger}
If $\phi \in \Out(F_r)$ is parageometric then $\lambda(\phi) > \lambda(\phi^\inv)$.
\end{theorem}

\begin{corollary}
\label{CorollaryParaInverse}
If $\phi \in \Out(F_r)$ is parageometric then $\phi^\inv$ is neither geometric nor parageometric.
\end{corollary}

\begin{proof} If $\phi^\inv$ is geometric then $\lambda(\phi^\inv) =
\lambda(\phi)$, whereas if $\phi^\inv$ is parageometric $\lambda(\phi^\inv) >
\lambda((\phi^\inv)^\inv) = \lambda(\phi)$.
\end{proof}

\begin{corollary}
\label{CorollaryTwoTrees}
A fully irreducible $\phi \in \Out(F_r)$ is geometric if and only if the $F_r$-trees $T_-$
and $T_+$ are both geometric.
\end{corollary}

\begin{proof} If $\phi$ is not geometric then $\phi^\inv$ is also not geometric, but if $T_-,T_+$ were
both geometric trees then it would follow by definition that $\phi$ and $\phi^\inv$ are both
parageometric, contradicting Corollary~\ref{CorollaryParaInverse}. The other direction was noted
above. 
\end{proof}

Corollaries~\ref{CorollaryParaInverse} and~\ref{CorollaryTwoTrees} have been proved independently
by Guirardel \cite{Guirardel:core}, by different means.

Theorem~\ref{TheoremLambdaBigger} and Corollary~\ref{CorollaryParaInverse} were first presented at
various seminars in the Fall of 2003, including the Topology Seminar at Princeton University. Our thanks
go to Baris Coskunuzer of that seminar for a question which quickly inspired the proof of the ``if''
direction of Corollary~\ref{CorollaryTwoTrees}.

\paragraph{Sketch of the proof of Theorem~\ref{TheoremLambdaBigger}.} After some preliminaries in Section~\ref{SectionPreliminaries}, in Section~\ref{SectionCharacterizing} we recall results from \OuterLimits\ that characterize when $T_+$ is a geometric $F_r$-tree: this happens if and only if some positive power of $\phi$ has a train track representative $g \from G \to G$ such that $G$ has a unique illegal turn,  $g$ has a unique periodic Nielsen path~$\rho$ (up to reversal), and $\Length(\rho) = 2 \Length(G)$. When $\rho$ exists, it is necessarily a \emph{fixed} Nielsen path, meaning that $\rho$ is fixed up to homotopy rel endpoints by the action of $g$. Also, $\rho$ decomposes at its illegal turn into $\rho = \alpha * \bar\beta$ where $\alpha,\beta$ are legal paths of length equal to~$\Length(G)$. In this situation, following \OuterLimits\ we construct a \nb{2}dimensional dynamical system $k\from K \to K$ representing $\phi$, which we call the \emph{wedge model}; a detailed description of the wedge model is given in Section~\ref{SectionWedgeModel}. The \nb{2}complex $K$ is obtained by attaching to $G$ a wedge $W$, a triangle with one side attached along $\alpha$ and the other side attached along~$\beta$. The unattached side of $W$ is vertical, and each vertical segment of $W$ has endpoints on a corresponding pair of points, one in $\alpha$ and one in $\beta$. The effect of $g \from G \to G$ is to fold $\rho$ by some amount, and this extends to a homotopy equivalence $k \from K \to K$ whose effect on $W$ is to collapse some vertical segments of $W$. The vertical segments of $W$ form leaf segments of a measured foliation $\fol^s$ on~$K$ called the \emph{stable foliation} of the wedge model $k \from K \to K$. The graph $G$ is transverse to the stable foliation, and the restriction to $G$ of the transverse measure on the stable foliation equals the Lebesgue measure along the train track $G$. Following \OuterLimits, in this situation we show that the attracting tree $T_+$ is the dual tree of the measured foliation $\fol^s$. Combining this construction with results of \BH, we show that $\phi$ is parageometric if and only if the two endpoints of $\rho$ are distinct. Moreover, in this case there exists an edge of $G$ which is covered exactly once by the Nielsen path $\rho$, and so is a free edge of the \nb{2}complex $K$. This simple observation, proved in Fact~\ref{FactParaWedgeModel}, plays a key role in the proof of Theorem~\ref{TheoremLambdaBigger}.

In Section~\ref{SectionStableFoliation} we carry out a detailed study of $\fol^s$, the stable foliation of the wedge model. We shall show that leaves of $\fol^s$ are trees, but they turn out to be trees of a rather thorny variety: these leaves have lots of valence~1 vertices, occuring in the interior of edges of $G$ that are free edges of the \nb{2}complex $K$. We are particularly interested in the collection of bi-infinite lines contained in leaves of $\fol^s$, which we denote $\Hull(\fol^s)$, the \emph{hull} of $\fol^s$. We shall use properties of $k \from K \to K$ to essentially identify $\Hull(\fol^s)$ with the expanding lamination of $\phi^\inv$.

Theorem~\ref{TheoremLambdaBigger} is proved in Section~\ref{SectionMainProof} by studying the asymptotic compression rate of $k$ on lines of $\Hull(\fol^s)$, defined to be the exponential growth rate in $n$ of the size of a subarc of a line of $\Hull(\fol^s)$ that $k^n$ collapses to a point. We compute this rate in several different ways. On the one hand, since $\Hull(\fol^s)$ is identified with the expanding lamination of $\phi^\inv$, the asymptotic compression rate of $k$ on lines of $\Hull(\fol^s)$ is equal to $\lambda(\phi^\inv)$. On the other hand, the abundance of valence~1 vertices in leaves of $\fol^s$ shows that the asymptotic compression rate of $k$ on lines in $\Hull(\fol^s)$ is strictly less than the exponential growth rate in $n$ for the size of a subtree of a leaf of $\fol^s$ that $k^n$ collapses to a point. The latter rate is simply the Perron-Frobenius eigenvalue of the transition matrix for $g \from G \to G$, which equals $\lambda(\phi)$. This is the culminating argument of the proof that $\lambda(\phi^\inv) < \lambda(\phi)$. 

\paragraph{When does $\lambda(\phi) = \lambda(\phi^\inv)$?} Consider a fully irreducible $\phi\in \Out(F_r)$. In the wake of our results one might wonder whether $\lambda(\phi)=\lambda(\phi^\inv)$ implies that $\phi$ is geometric. Here is an easy construction of counterexamples: fully irreducible outer automorphisms $\phi$ which are not geometric and yet which satisfy $\lambda(\phi)=\lambda(\phi^\inv)$; by Theorem~\ref{TheoremLambdaBigger}, neither is such a $\phi$ parageometric. 

In any group~$G$, if $g,g' \in G$ have order~2 then $gg'$ and $(gg')^\inv$ are conjugate. It follows that if $\Psi,\Psi' \in \Out(F_r)$ have order~2, and if $\phi=\Psi\Psi'$ is fully irreducible, then $\lambda(\phi) =
\lambda(\phi^\inv)$. 

For a concrete example, let $\Psi \in \Out(F_3)$ be the order two element represented by the automorphism $a \mapsto b$, $b \mapsto a$, $c \mapsto c$. Let $\Psi' = \Theta\Psi\Theta^\inv$ be a conjugate of~$\Psi$. If the conjugating element $\Theta$ is picked randomly then one might expect that $\phi=\Psi\Psi'$ is fully irreducible and has a train track representative with no periodic Nielsen paths, and so $\phi$ is  nongeometric. Taking $\Theta$ to be the fourth power of the outer automorphism $a \mapsto b$, $b \mapsto c$, $c \mapsto \bar b a$ considered earlier, and applying the train track algorithm of \BH, one obtains the following train track map $g \from G \to G$ representing $\phi$. The graph $G$ has two vertices $r,q$, four edges $B,C,D,E$, with $B$ from $q$ to $r$, $C$ from $r$ to $q$, $D$ from $r$ to $q$, and $E$ from $q$ to $q$, and $g$ is defined by $g(B) = C \bar E \bar C D E$, $g(C) = \bar C \bar B \bar E \bar D$, $g(D)=B$, and $g(E)=CB$. The expansion factor is $\lambda = 3.199158087\ldots$. 

To verify that $g$ represents a nongeometric, fully irreducible outer automorphism it is sufficient to check three things. First,
the transition matrix of $g$ is positive. Second, at each of the two vertices $v=r,q$ of $G$, the graph of turns taken at $v$ is connected; this is the graph with one vertex for each oriented edge with initial vertex $v$, and one edge for each pair of oriented edges $E_1 \ne E_2$ with initial vertex $v$ such that the $g$ image of some oriented edge of $G$ contains the subpath $\bar E_1 E_2$. Third, $g$ has no periodic Nielsen paths, which can be checked by the following expedient. Factor $g$ into Stallings folds, $G=G_0 \to G_1 \to \cdots \to G_n=G$; we did this with $n=8$. Let $\E_0$ be the set of length 2 edge paths in $G_0$ with an illegal turn; there are two such paths up to reversal, $E \bar D$ and $\bar E B$. For negative integers $i$ define $\E_i$ inductively to be a set of edge paths in $G_i$ (with indices taken modulo $n$) each with one illegal turn, as follows: for each $\gamma\in \E_i$, take all paths in $G_{i-1}$ with exactly one illegal turn whose straightened image in $G_i$ is $\gamma$, and put each such path in $\E_{i-1}$. Carrying this process out, we computed that the set $\E_{-12}$ is empty. This shows that $g$ has no Nielsen paths, and so it represents a fully irreducible, nongeometric outer automorphism. 

Just as a check, we also inverted the sequence of Stallings folds and applied the train track algorithm to verify that the expansion factor of the inverse is also equal to $3.199158087\ldots$.

To get wider classes of examples, consider a fully irreducible $\phi \in\Out(F_r)$ with expanding lamination denoted $\Lambda^u$. By Section~2 of \BookZero\, the group $\Out(F_r)$ acts on the set of expanding laminations of fully irreducible elements, and there is a homomorphism $\ell^u\from\Stab(\Lambda^u) \to\reals_+$ with discrete image and finite kernel such that
$\ell^u(\phi)=\lambda(\phi)$, and $\ell^u(\Psi) = \lambda(\Psi)$ as long as $\ell^u(\Psi)\ge 1$.  Applying this to $\phi^\inv$ with expanding lamination denoted $\Lambda^s$, we obtain a homomorphism $\ell^s\from \Stab(\Lambda^s) \to\reals_+$ such that $\ell^s(\phi^\inv) = \lambda(\phi^\inv)$ and $\ell^s(\Psi) = \lambda(\Psi)$ as long as $\ell^s(\Psi)\ge 1$. Applying Proposition~2.16 of \BookOne\ it follows that $\Stab(\Lambda^s)=\Stab(\Lambda^u)$, a subgroup of $\Out(F_r)$ that we denote $V_\phi$. From the properties of the homomorphisms $\ell^u,\ell^s \from V_\phi \to\reals_+$ it follows that the infinite cyclic group $\<\phi\>$ has finite index in $V_\phi$, and so any two elements of $V_\phi$ not contained in the common kernel of $\ell^s,\ell^u$ have nonzero powers that are equal. This implies that if $\lambda(\phi)=\lambda(\phi^\inv)$ then $\lambda(\Psi)=\lambda(\Psi^\inv)$ for any $\Psi\in V_\phi$. Note that $V_\phi$ is the \emph{virtual centralizer} of $\<\phi\>$ in $\Out(F_r)$, consisting of all $\Psi\in\Out(F_r)$ that commute with some positive power of $\phi$; $V_\phi$ is contained in the virtual centralizer because $\<\phi\>$ has finite index in $V_\phi$; and if $\Psi\not\in V_\phi$ then $\Psi\phi^k\Psi^\inv \ne \phi^k$ because their attracting fixed points in $\bdy\X_r$ are distinct, by Proposition~2.16 of \BookOne. 

It might be interesting to find necessary and sufficient conditions for the condition $\lambda(\phi) =
\lambda(\phi^\inv)$, for fully irreducible $\phi \in \Out(F_r)$. For example, is it necessary that
$V_\phi$ contains a fully irreducible element that is either geometric or conjugate to its own
inverse? 

At the very least, it would seem that the property $\lambda(\phi)\ne\lambda(\phi^\inv)$ is generic, and the property that $\phi$ be fully irreducible with nongeometric fixed trees $T_-$, $T_+$ is also generic. We invite the reader to take a random word of, say, twenty or more Nielsen generators of $\Out(F_3)$ and verify that the resulting outer automorphism $\phi$ is fully irreducible, neither geometric nor parageometric, and satisfies $\lambda(\phi) \ne \lambda(\phi^\inv)$.

\section{Preliminaries}
\label{SectionPreliminaries}

\subsection{Outer automorphisms and outer space.}

The definitions in this section follow several sources. For the foundations of marked graphs,
$\reals$-trees, and outer space, including many of the facts recalled below without citation, see
\cite{CullerMorgan:Rtrees} and \cite{CullerVogtmann:moduli}. For concepts of irreducibility see \BH. A
good overview is given in \cite{Vogtmann:OuterSpaceSurvey}.

\paragraph{Outer automorphisms of free groups.} Fix an integer $r \ge 2$, let $F_r$ denote the free group
of rank~$r$, let $\Out(F_r) = \Aut(F_r) / \Inn(F_r)$ denote its outer automorphism group, and let
$\C$ denote its set of nontrivial conjugacy classes. Let $R_r$ denote the rose
with $r$-petals and identify $\pi_1(R_r) \approx F_r$, so the group $\Out(F_r)$ is identified with the
group of homotopy classes of self-homotopy equivalences of $R_r$. Given $\phi \in \Out(F_r)$ let $f_\phi
\from R_r \to R_r$ be a representative homotopy equivalence. $\Out(F_r)$ acts naturally on $\C$, and on
conjugacy classes of subgroups of $F_r$. We say that $\phi \in \Out(F_r)$ is \emph{reducible} if there
exists a nontrivial free factorization $F_r = A_1 * \cdots * A_k * B$ so that $\phi$ permutes the
conjugacy classes of $A_1,\ldots,A_k$. If $\phi$ is not reducible then it is \emph{irreducible}. If
$\phi^n$ is irreducible for all $n \ge 1$ then $\phi$ is \emph{fully irreducible}\footnote{called
``irreducible with irreducible powers'' or ``IWIP'' in the literature}. Note that $\phi$ is irreducible
if and only if $\phi^\inv$ is, and the same for complete irreducibility.

\paragraph{Outer space and its boundary.} An \emph{$F_n$-tree} is an $\reals$-tree $T$ equipped with an
action of $F_n$ that is minimal (no proper nonempty subtree is invariant) and nonelementary ($T$ is not
a point or a line). An $F_r$-tree is \emph{proper} if the action is properly discontinuous, and it is
\emph{simplicial} if $T$ is a simplicial complex.  Two $F_n$ trees are \emph{isometrically (resp.
homothetically) conjugate} if there is an isometry (resp. homothety) between them that conjugates one
action to the other. \emph{Outer space} $\X_r$ is the set of homothetic conjugacy classes of proper,
simplicial $F_n$-trees, with topology induced by embedding $\X_r \to \P\reals^\C$ as follows: first
embed the set of isometric conjugacy classes into $\reals^\C$ using translation length as a class
function on $F_n$, and then projectivize. The image of this embedding is precompact, and its closure and
boundary are denoted $\overline\X_r$ and $\bdy\X_r = \overline\X_r - \X_r$. 

Points of outer space can also be represented as marked graphs, as follows. A marked graph is a graph
$G$ with all vertices of valence $\ge 3$, equipped with a path metric, and with a homotopy equivalence
$R_r \to G$ called the \emph{marking}. If a base point $p \in G$ happens to be imposed, the homotopy
class of a marking of $G$ determines and is determined by an isomorphism $F_r\to\pi_1(G,p)$ up to
precomposition by an inner automorphism of $F_r$. Two marked graphs $G,G'$ are \emph{isometric (resp.
homothetic)} if there exists an isometry (resp. homothety) $G\to G'$ which, together with the markings,
makes the following diagram commute up to homotopy:
$$\xymatrix{G\ar@/^1pc/[rr] & R_r \ar[l] \ar[r] & G'}
$$
Passage to the universal covering space induces a bijection between the set of homothety classes of
marked graphs and the set $\X_r$. The embedding $\X_r \to \P\reals^\C$ can be understood by first
associating to a marked graph $G$ the class function on $F_n$ that associates to an element of
$F_n$ the length of the shortest loop in $G$ representing the free homotopy class of that element, and
then projectivizing.

The length of an object in a geodesic metric space is denoted $\Length(\cdot)$, with a subscript to denote the metric space when the context is not clear, for example $\Length_G(\cdot)$ in the marked graph $G$. Also, when a marked graph $G$ is clear from the context then we use the constant $L$ to denote $\Length(G)$.

The group $\Out(F_r)$ acts on $\overline\X_r$ on the right, as follows. Let $[\cdot]$ denote
homothety class. For each $F_r$-tree $T$ one can precompose the action $F_r\to \Isom(T)$ with an
automorphism $F_r \to F_r$ representing $\phi$, to get $[T] \phi$. In terms of a marked graph $G$, one
can precompose the marking $R_r \to G$ with a homotopy equivalence $R_r \to R_r$ that represents $\phi$,
to obtain $[G]\phi$. This action preserves the topology, and it preserves outer space $\X_r$ itself and
its boundary $\bdy \X_r$.

\paragraph{Source--sink dynamics.} If $\phi \in \Out(F_r)$ is fully irreducible then there exist
$T_- \ne T_+ \in \bdy\X_r$ such that for every $x\in \overline\X_r$,
$$\lim_{n \to -\infinity} \phi^n(x) = T_- \qquad\text{and}\qquad \lim_{n
\to +\infinity}\phi^n(x) = T_+
$$
where these limits take place in $\overline\X_r$. This was proved for $x \in \X_r$ in
\cite{BFH:laminations} and extended to all $x \in \overline\X_r$ in \cite{LevittLustig:NorthSouth}. We
call $T_-$ the \emph{repelling tree} and $T_+$ the \emph{attracting tree} of~$\phi$.

\paragraph{Geometric $F_r$-trees.} We review measured foliations on \nb{2}complexes
and geometric trees following \cite{LevittPaulin:geometric}.

Consider a connected simplicial \nb{2}complex $K$ which is not a point. A \emph{measured foliation} on
$K$ is defined by specifying measured foliations on each \nb{2}simplex of $K$ which fit together
compatibly along \nb{1}simplices. To be precise, a measured foliation on a~1 or~2 dimensional simplex
$\sigma$ in $K$ is determined by choosing a simplicial homeomorphism $\sigma' \to \sigma$ where
$\sigma'$ is a rectilinear simplex in $\reals^2$, and pushing forward the vertical foliation on
$\reals^2$ with the transverse measure $\abs{dx}$. A \emph{leaf segment} in $\sigma$ is the pushforward
of $\sigma'$ intersected with a vertical line; for example, if $\sigma$ is a \nb{1}simplex then either
$\sigma$ is a single leaf segment or each point of $\sigma$ is a leaf segment. A measured foliation on a
\nb{2}simplex restricts to a measured foliation on each of its edges. A measured foliation on $K$ is
determined by choosing a measured foliation on each \nb{2}simplex of $K$, so that for each \nb{1}simplex
$e$, all of the measured foliations on $e$ obtained by restricting to $e$ from a \nb{2}simplex incident
to $e$ agree with each other.

We will often suppress the simplicial structure on $K$, so a measured foliation on a cell
complex means, formally, a measured foliation on some simplicial subdivision.

Let $\fol$ denote a measured foliation on $K$, and $\fol\restrict\sigma$ its restriction to each
simplex $\sigma$ of $K$. The collection of leaf segments in \nb{1}simplices and \nb{2}simplices define
a relation on $K$, two points being related if they are contained in the same leaf segment. This
relation generates an equivalence relation on $K$. The equivalence classes are called \emph{leaves}.
The leaf containing a point $x \in K$ can be built up inductively as follows: let $\ell_1$ be the union
of all leaf segments containing $x$; for $i \ge 1$ let $\ell_{i+1}$ be the union of all leaf segments
containing points of $\ell_i$; finally, the leaf containing $x$ is $\union_{i=1}^\infinity \ell_i$.  

Given a measured foliation $\fol$ on $K$ and a path $\gamma \from I \to K$, pulling back the
tranverse measure locally gives a measure on $I$, whose integral is denoted $\int_\gamma\fol$. 

Given a finite \nb{2}complex $K$ with measured foliation $\fol$ and a surjective homomorphism $h \from
\pi_1(K)\to F_r$, let $\wt K \to K$ be the covering space corresponding to $\ker(h)$, and let
$\wt\fol$ be the lifted measured foliation on $\wt K$.  Define a pseudo-metric on $\wt K$ where
$d(x,y)$ is the infimum of the transverse measures of paths from $x$ to~$y$. Let $T$ be the
associated metric space, whose points are the equivalence classes determined by the relation
$d(x,y)=0$. Note that if $x,y$ are in the same leaf of $\wt K$ then $d(x,y)=0$, but the converse need
not hold in general. The action of $F_r$ on $\wt K$ induces an isometric action of~$F_r$ on the
metric space $T$. We assume that each \nb{1}simplex $e$ of $\wt K$ is a geodesic, that is, if
$\bdy e = \{v,w\}$ then $d(v,w) = \int_e \wt\fol$. Under this assumption, Levitt and Paulin
\cite{LevittPaulin:geometric} prove that $T$ is an $F_r$-tree, called the \emph{dual $F_r$-tree}
of the measured foliation~$\fol$.

An $F_r$-tree $T$ is said to be \emph{geometric} if there exists a finite \nb{2}complex $K$ with
measured foliation $\fol$, and a surjective homomorphism $\pi_1(K)\to F_r$, such that each
edge of $\wt K$ is a geodesic, and such that $T$ is isometrically conjugate to the dual $F_r$-tree
of $\fol$. 

\paragraph{Geometric and parageometric outer automorphisms.} An outer automorphism $\phi \in \Out(F_r)$
is \emph{geometric} if there exists a compact surface $S$, an isomorphism $F_r \approx \pi_1 S$, and a
homeomorphism $h \from S \to S$, such that the outer automorphism of $F_r$ induced by $h$ is equal to
$\phi$. If $\phi$ is fully irreducible and geometric, then its attracting and repelling $F_r$-trees
$T_+,T_-$ are both geometric $F_r$-trees.

Consider now a fully irreducible $\phi \in \Out(F_r)$ with attracting $F_r$-tree $T_+ \in
\bdy\X_r$. We say that $\phi$ is \emph{parageometric} if $T_+$ is a geometric $F_r$-tree but $\phi$
is not a geometric outer automorphism.

\subsection{Train tracks and laminations}

The definitions in this section follow \BH\ and \BookZero.

\paragraph{Topological representatives and Markov partitions.} Given $\phi \in \Out(F_r)$, a marked
graph $G$ with marking $\mu \from R_r \to G$, and a homotopy equivalence $g \from G \to G$, we say that
$g$ is a \emph{topological representative} of $\phi$ if $g$ takes vertices to vertices, $g$ is an
immersion on each edge, and the following diagram commutes up to homotopy:
$$\xymatrix{
R_r \ar[r]^{f_\phi} \ar[d]_\mu & R_r \ar[d]^\mu \\
G \ar[r]_g & G
}
$$
and so the composition $R_r \to G \xrightarrow{g} G$ represents the point $[G]\phi\in \X_r$. The set of
edges $\Edges$ of $G$ forms a Markov partition for $g$, meaning that for any $e,e' \in \E$, each
component of $e\intersect g^\inv(\interior(e'))$ is mapped by $g$ homeomorphically onto
$\interior(e')$. The \emph{transition graph} $\TG$ of $g$ is a directed graph whose vertex set is the
set $\Edges$, such that for each $e,e' \in \Edges$, the directed edges from $e$ to $e'$ are in
one-to-one correspondence with the components of $e\intersect g^\inv(\interior(e'))$. The
\emph{transition matrix} of $g$ is the function $M\from\Edges \cross \Edges\to \Z$, where
$M(e,e')$ equals the number of directed edges in $\TG$ from $e$ to $e'$, in other words, the number
of times that $g(e')$ traverses $e$ in either direction. Note that $M^n(e,e')$ is the number of directed
paths from $e$ to $e'$ of length $n$. We say that $\TG$, $M$, and $g \from G \to G$ are
\emph{irreducible} if there is an oriented path from any vertex of $\TG$ to any other vertex,
equivalently, for all $e,e'$ there exists $n$ such that $M^n(e,e') \ne 0$. If irreducibility holds then
the Perron-Frobenius theorem implies that there exists a unique $\lambda \ge 1$ such that $M$ has a
positive (right) eigenvector with eigenvalue $\lambda$. If $M$ has the stronger property that some
positive power has all positive entries then $\lambda>1$ and a positive eigenvector is unique up to
positive scalar multiple.

\paragraph{Train tracks.} A \emph{direction} of $G$ at a vertex $v$ is the germ, up to
reparameterization, of an immersed path with initial point $v$. Each direction is uniquely represented
by an oriented edge $e$ with initial point $v$, but we occasionally use other paths with initial point
$v$ to represent directions. A \emph{turn} of $G$ at $v$ is an unordered pair
$\{e,e'\}$ of directions at $v$; the turn is \emph{nondegenerate} if $e \ne e'$, otherwise the turn is
\emph{degenerate}. 

An \emph{edge path} in $G$ will always mean a concatenation of the form $\gamma = e_0 * e_1 * \cdots *
e_{k-1} * e_k$, $k \ge 0$, where $e_1,\ldots,e_{k-1}$ are oriented edges and $e_0,e_k$ are subsegments of
oriented edges. Often we say ``path'' when ``edge path'' is meant; the context should make this clear.
Given an edge path $\gamma \from I \to G$ and $t\in\interior(I)$ so that $\gamma(t)$ is a vertex of $G$,
let $e,e'$ be the two directions of $\gamma$ at this point, that is: subdivide at $t$ to obtain a
concatenation $\gamma = \alpha * \beta$, let $e$ be the direction of $\bar\alpha$ at its initial point,
and let $e'$ be the direction of $\beta$ at its initial point. With this notation we say that $\gamma$
\emph{takes the turn $\{e,e'\}$ at the parameter value $t$}. If $t$ is understood then we just say that
$\gamma$ \emph{takes the turn $\{e,e'\}$}. 

A topological representative $g \from G \to G$ acts on the set of directions and on the set of turns of
$G$. A nondegenerate turn is \emph{illegal} if its image under some positive power of $g$ is
degenerate, otherwise the turn is \emph{legal}. An edge path $\alpha$ is \emph{legal} if every turn taken
by $\alpha$ is legal, in particular every legal path is immersed. Given a path $\alpha$ in $G$, let
$\alpha_\#$ denote the immersed path (or constant path) which is homotopic to $\alpha$ rel endpoints, so
if $\alpha$ is immersed then $\alpha_\#=\alpha$, and if $\alpha$ is legal then
$g^n(\alpha)_\#=g^n(\alpha)$ for all $n\ge 0$. 

A topological representative $g \from G \to G$ of $\phi \in \Out(F_r)$ is a \emph{train track
representative} of $\phi$, and $g$ is a \emph{train track map}, if for each edge $E$ of $G$,
the map $g \restrict E \from E \to G$ is a legal path, equivalently, $g^n \restrict E$ is an immersion
for each $n\ge 1$. A train track map $g$ is \emph{fully irreducible} if $g^n$ is irreducible for
all $n \ge 1$.

\begin{proposition}[\BH] 
\label{PropStretchFactor}
If $\phi \in \Out(F_r)$ is fully irreducible then $\phi$ has a fully irreducible train track
representative $g \from G \to G$. The transition matrix $M_{g}$ has a positive power, and so there
exists $\lambda(g) > 1$ and vector $v \from\Edges \to \reals$, with $\lambda(g)$ unique and $v$
unique up to a positive scalar multiple, so that $Mv = \lambda(g) v$. We also have
$\lambda(\phi)=\lambda(g)$ (see Remark~1.8 of \BH).
\end{proposition}

With $g \from G \to G$ as in this proposition, we may assign a path metric to $G$, also called the Lebesgue measure on $G$, so that each edge $e$ has length $\Length_G(e)=v(e)$. We may then alter $g$ on each edge $e$ by a homotopy rel endpoints so that $g\restrict e$ stretches path length by a constant factor of $\lambda(g)$; the resulting map is still a train track map. The number $\lambda(g)=\lambda(\phi)$ is called the \emph{stretch factor} of $g$.

Henceforth we \emph{always assume} without comment that if $g \from G \to G$ is a train track
representative of a fully irreducible $\phi \in \Out(F_r)$ then $g$ is fully irreducible and
$g$ stretches path length on $G$ by the constant factor $\lambda(\phi)$. It follows that for any legal path $\gamma$ the legal path $g(\gamma)_\# = g(\gamma)$ has length equal to $\lambda(\phi) \, \Length(\gamma)$.

\paragraph{The geodesic lamination of a free group.} In this heading and the next we review the results
from \BookZero\ concerning the construction and properties of the \emph{expanding} or \emph{unstable}
lamination $\Lambda^u(\phi)$ of a fully irreducible $\phi \in \Out(F_r)$. We use here a slightly
different point of view than in \BookZero, presenting laminations as Hausdorff objects rather than
non-Hausdorff.

Consider a marked graph $G$ with universal cover $T$. The \emph{geodesic lamination} of $T$, denoted
$\Lambda_T$, is the set of pairs $(\ell,x)$ where $\ell \subset T$ is a bi-infinite, unoriented line and
$x\in\ell$, equipped with the compact open topology where a neighborhood $U_\epsilon$ of $(\ell,x)$ is
the set of all $(\ell',x')$ such that $d(x,x') \le \epsilon$ and $\ell\intersect \ell'$ contains a
$1/\epsilon$ neighborhood of $x$ in $\ell$ and a $1/\epsilon$ neighborhood of $x'$ in $\ell'$. The
\emph{projection map} $\Lambda_T \to T$ is the map $(\ell,x) \to x$. A \emph{leaf} of $\Lambda_T$
corresponding to a bi-infinite line $\ell$ is the set of all $(\ell,x)$ such that $x \in \ell$; we
shall often confuse a bi-infinite line in $T$ with its corresponding leaf. A point of $\Lambda_T$
can also be described as a geodesic embedding $\gamma \from \reals \to T$ modulo precomposition by the
involution $x\leftrightarrow -x$, where the point of $\Lambda_T$ corresponding to $\gamma$ is the pair
$(\gamma(\reals),\gamma(0))$. The action of $F_r$ on $T$ induces a properly discontinuous, cocompact
action on~$\Lambda_T$ whose quotient space, a compact lamination denoted $\Lambda_G$, is the
\emph{geodesic lamination} of $G$. The projection map $\Lambda_T \to T$ descends to a projection map
$\Lambda_G \to G$. An element of $\Lambda_G$ can also be described as a locally geodesic immersion
$\reals \to G$ modulo precomposition by the involution $x \leftrightarrow -x$ on the parameter
domain~$\reals$. A \emph{sublamination} of any lamination is a closed subset that is a union of leaves.
By compactness of $\Lambda_G$, every sublamination of $\Lambda_G$ is compact. 

Given two marked graphs $G,G'$, any homotopy equivalence $g \from G \to G'$ induces a homeomorphism
$g_* \from \Lambda_G \to \Lambda_{G'}$ well defined up to isotopy, defined as follows. First alter $g$
by homotopy so that it takes vertices to vertices and is affine on each edge. Lift $g$ to the
universal covers $\tilde g \from T = \wt G \to \wt G' = T'$. There is an automorphism $\Phi \from F_r
\to F_r$ such that $g$ satisfies $\Phi$-twisted equivariance, meaning that $f(g(x)) = g(\Phi(f)(x))$
for all $x\in T$, $f \in F_r$. Consider a leaf $\ell$ of $\Lambda_T$. Since $\tilde g$ is a
quasi-isometry, $\tilde g(\ell)$ is a quasi-geodesic embedding of $\reals$, and so the image $\tilde
g(\ell)$ has finite Hausdorff distance from some leaf that we shall denote $\tilde g_*(\ell)$. Define
a function $\tilde g_\# \from \Lambda_T \to \Lambda_{T'}$ that maps each leaf $\ell$ to $\tilde
g_*(\ell)$, by postcomposing the map $\tilde g$ with the closest point projection from $\tilde
g(\ell)$ onto $\tilde g_*(\ell)$. The map $\tilde g_\#$ is continuous, $\Phi$-twisted equivariant, and
induces a bijection of leaves. The image of the map $\ell \to \tilde g(\ell)$ is the line $\tilde
g_*(\ell)$ union a disjoint set of finite trees attached to the line, and the effect of the closest
point projection is to collapse each of these finite trees to the point where it attaches to the line;
it follows that the map $\tilde g_\#$ is \emph{leafwise monotonic}, meaning that for each leaf $\ell$
the map $\ell \xrightarrow{\tilde g_\#} \tilde g_*(\ell)$ has the property that each point pre-image is
an arc. We can now perturb $\tilde g_\#$ to get a homeomorphism $\tilde g_*$, still satisfying
$\Phi$-twisted equivariance, and $\tilde g_*$ is well-defined up to $\Phi$-twisted equivariant isotopy.
It follows that $\tilde g_*$ descends to the desired homeomorphism $g_*$, well-defined up to isotopy. 

For any $g \from G \to G'$ as above that preserves the marking (in the sense that the marking $R_r \to
G$, postcomposed with $g$, is homotopic to the marking $R_r \to G'$), note that the map $\tilde g_*
\from \Lambda_G \to \Lambda_{G'}$ is natural in the sense that for any two homotopy equivalences $G
\xrightarrow{g} G'\xrightarrow{g'} G''$ that preserve markings, the composition $\Lambda_G
\xrightarrow{g_*} \Lambda_{G'} \xrightarrow{g'_*} \Lambda_{G''}$ is isotopic to $(g' \composed g)_*$.
We are therefore justified in talking about ``the'' geodesic lamination $\Lambda_r$ of $F_r$, as
represented by $\Lambda_G$ for any marked graph $G$. We are also justified in talking about
a sublamination $\Lambda' \subset \Lambda_r$ of the geodesic lamination of $F_r$, represented as a
sublamination $\Lambda'_G \subset \Lambda_G$ for any marked graph $G$, with the property that for any
marking preserving homotopy equivalence $g \from G \to G'$ we have $g_*(\Lambda'_G) = \Lambda'_{G'}$.
We also say that $\Lambda'_G$ is the \emph{realization} of $\Lambda'$ in the marked graph $G$.

\paragraph{Expanding laminations.} Let $\phi \in \Out(F_r)$ be fully irreducible.
The \emph{expanding} or \emph{unstable lamination} of $\phi$ is a sublamination $\Lambda^u(\phi)$ of
the geodesic lamination of $F_r$, defined as follows. Choose any train track representative $g\from G
\to G$. Choose an edge $e$ and a periodic point $x\in\interior(e)$ of periodicity $p$. As $n\to
\infinity$, the maps $e\xrightarrow{g^{np}} G$ can be reparameterized as a nested sequence of isometric
immersions of larger and larger subintervals of $\reals$, each interval containing $0$ and each immersion
taking $0$ to $x$. The union of these immersions is a bi-infinite geodesic in $G$, that is, a leaf of
$\Lambda_G$. The closure of this leaf in $\Lambda_G$ is defined to be the realization of
$\Lambda^u=\Lambda^u(\phi)$ in the marked graph $G$, denoted $\Lambda^u_G$.
$\Lambda^u$ is well-defined, independent of the choice of $e$ and $x$, and also independent of $g$
meaning that for any other train track representative $g' \from G' \to G'$ and any marking preserving
homotopy equivalence $h\from G\to G'$ we have $h_*(\Lambda^u_G) = \Lambda^u_{G'}$. 

The expanding lamination $\Lambda^u(\phi)$ is \emph{minimal}, meaning that its only nonempty
sublamination is itself, in other words, every leaf is dense; see Section~1 of \BookZero. Also, the projection from $\Lambda^u_G$ to $G$ is surjective for any train track representative $g \from G \to G$, because some power of $g$ has positive transition matrix.

Note that the action of $g_*$ on $\Lambda_G$ restricts to an action on $\Lambda^u_G$ which expands length
by the exact factor of $\lambda(\phi)$, that is, $\Length(g_*(\ell)) = \lambda(\phi) \Length(\ell)$ for
any leaf segment $\ell$ of $\Lambda^u_G$.

\section{Geometric trees and the wedge model}
\label{SectionCharacterizing}

Let $\phi \in \Out(F_r)$ be fully irreducible with attracting $F_r$-tree $T_+$. The primary goal of this section is to characterize when $T_+$ is a geometric tree, and to use this characterization as an opportunity for introducing the wedge model of $\phi$, which will play such an important role in later sections. We will also obtain a characterization of parageometricity of $\phi$. Both characterizations are stated in terms of train track representatives of $\phi$ and its positive powers. These characterizations are restatements and reworkings of results in \OuterLimits\ and in \BH. We provide full details of proof, in part because of the limited availability of \OuterLimits, but also because of our need to develop a complete description of the wedge model.

The reader who wants to skip quickly to the definition of the wedge model should first read Section~\ref{SectionNielsen} and then Section~\ref{SectionWedgeModel}, skipping Sections~\ref{SectionDirectLimits} and~\ref{SectionNoNielsen}, although Section~\ref{SectionDirectLimits} will be needed to understand geometricity of $T_+$.

\subsection{Nielsen paths} 
\label{SectionNielsen}

Consider a train track map $g \from G \to G$. A \emph{(fixed) Nielsen path} of $g$ is a locally geodesic path $\rho\from [a,b] \to G$ such that $\rho(a)$ and $\rho(b)$ are fixed points, and $g \composed \rho$ is homotopic rel endpoints to $\rho$. A \emph{periodic Nielsen path} of $g$ is a Nielsen path of some power $g^n$ with $n \ge 1$. A (periodic) Nielsen path $\rho$ is \emph{indivisible} if it cannot be written as a nontrivial concatenation of (periodic) Nielsen paths. Assuming that $g$ is irreducible, every indivisible periodic Nielsen path has a unique \emph{legal decomposition} $\rho = \alpha * \bar \beta$ where $\alpha,\beta \from [0,b] \to G$ are legal paths of equal length, and $\rho$ takes an illegal turn at the concatenation point. Depending on the context we may write this decomposition in other forms, for example, $\rho=\alpha_1 * \overline\alpha_2$. For details on Nielsen paths see \BH.

The following result is essentially proved in \BH\ using the theory of stable train tracks. We say that a surjective map $\alpha \from A \to B$ is \emph{generated by} a relation $\sigma \subset A \cross A$ if the equivalence relation generated by $\sigma$ has equivalence classes identical to the point pre-images of $f$.

\begin{proposition}
\label{PropNielsenUnique}
If $\phi \in \Out(F_r)$ is fully irreducible then there exists $k \ge 1$ and a train track representative $g \from G \to G$ of $\phi^k$ such that one of the following holds:
\begin{enumerate}
\item \label{ItemNoNielsen}
$g$ has no periodic Nielsen paths. 
\item \label{ItemNielsenUnique}
\begin{enumerate}
\item \label{ItemRho}
$g$ has an indivisible periodic Nielsen path $\rho$ unique up to orientation reversal, $\rho$ is a fixed Nielsen path, the image of $\rho$ is all of $G$, the illegal turn taken by $\rho$ is the unique illegal turn of $g$, and $\Length(\rho)=2\Length(G)$. 
\item \label{ItemGeneration}
Letting the legal decomposition be $\rho = \alpha * \bar\beta$ with $\alpha,\beta \from [0,\Length(G)] \to G$, the map $g$ is generated by the relation 
$$\alpha(t) \sim \beta(t)\quad\text{for}\quad t \in [\Length(G)/\lambda,\Length(G)]
$$
\end{enumerate}
\end{enumerate}
\end{proposition}

A train track map $g \from G \to G$ satisfying~(\ref{ItemNielsenUnique}) is said to be \emph{Nielsen unique}. 

\begin{proof} We briefly review the theory of stable train tracks. Let $g \from G \to G$ be a train track representative of $\phi$, let $\rho = \alpha * \bar \beta$ be an irreducible Nielsen path, and suppose that the illegal turn of $\rho$ is immediately folded by $g$. Now apply Stalling's method of factoring $g$: let $E_\alpha,E_\beta$ be maximal initial oriented segments of $\bar\alpha,\bar\beta$, respectively, such that each of $E_\alpha,E_\beta$ is contained in a single edge of $G$, and the paths $g(E_\alpha)$ and $g(E_\beta)$ in $G$ are the same; and then factor $g$ as $g=h \composed f$ where $G \xrightarrow{f} G' \xrightarrow{h} G$, the graph $G'$ is the quotient of $G$ obtained by identifying $E_\alpha$ and $E_\beta$ to a single arc, $f \from G \to G'$ is the quotient map, and $h \from G' \to G$ is induced by $g$ under the quotient map $f$. The map $g' = f \composed h \from G' \to G'$ is a train track representative of $\phi$, obtained from $g$ by \emph{folding the irreducible Nielsen path $\rho$}. This fold is said to be \emph{full} if either $E_\alpha$ or $E_\beta$ is an entire edge. Note that a sufficient condition for fullness to fail is if the Nielsen path $\rho = \alpha * \bar\beta$ is \emph{small} meaning that each of the paths $\alpha,\bar\beta$ is a subarc of some edge. We note that the indivisible periodic Nielsen paths of $g$ and of $g'$ are in one-to-one, periodic preserving correspondence: paths $\sigma,\sigma'$ correspond in this way if $\sigma' = (f \composed \sigma)_\#$.

A train track representative $g \from G \to G$ of $\phi$ is \emph{unstable} if there exists a sequence of train train representatives $g=g_0,g_1,\ldots,g_k$ such that each $g_i$ is obtained from the previous $g_{i-1}$ by folding an irreducible Nielsen path, and 
the fold from $g_{k-1}$ to $g_k$ is not full. If $g$ is not unstable then $g$ is \emph{stable}. It follows that any train track representative obtained from a stable train track representative by folding an irreducible Nielsen path is also stable. Stability of $g$ implies that one of the following statements holds:
\begin{itemize}
\item[(\ref{ItemNoNielsen}$'$)] $g$ has no indivisible fixed Nielsen paths.
\item[(\ref{ItemRho}$'$)] $g$ has an indivisible fixed Nielsen path $\rho = \bar\alpha * \beta$ unique up to orientation reversal, the image of $\rho$ is all of $G$, the illegal turn in $\rho$ is the unique illegal turn of $g$, and $\Length(\rho)=2\Length(G)$. 
\end{itemize}
For these statements see Section~3 of \BH, particularly Lemma~3.9 of \BH, and for the proof that the image of $\rho$ is all of $G$ see the top of page~28 of \BH. 

We turn to the proof of items~(\ref{ItemNoNielsen}) and~(\ref{ItemRho}). Let $g' \from G' \to G'$ be any train track representative of $\phi$. From \BookOne\ Lemma 4.2.5, $g'$ has only finitely many indivisible periodic Nielsen paths. Pass to a power of $g'$ so that all indivisible periodic Nielsen paths of $g'$ are fixed. In \BH, a procedure is described which, from $g'$, produces a stable train track representative $g \from G \to G$, and so $g$ satisfies items~(\ref{ItemNoNielsen}$'$) and~(\ref{ItemRho}$'$). But we need to replace the word ``fixed'' by the word ``periodic'' in these items. We do this by using some of the details of the stabilization procedure, described in Section~3 of \BH, which produces $g$ from $g'$. The output of the stabilization procedure is a sequence of train track representatives $g'=g_0,g_1,\ldots,g_N=g$ so that for each $n=1,\ldots,N$ one of two possibilities holds. In one case, $g_n$ is obtained from $g_{n-1}$ by folding an invidisible Nielsen path of $g_{n-1}$; in this case the indivisible periodic Nielsen paths of $g_{n-1}$ and of $g_n$ are in one-to-one, period preserving correspondence, as noted above. In the other case, $g_n$ is obtained from $g_{n-1}$ by an operation that eliminates a small indivisible Nielsen path, and the remaining indivisible periodic Nielsen paths of $g_{n-1}$ that are not eliminated are in one-to-one period preserving correspondence with the indivisible periodic Nielsen paths of $g_n$. Since all indivisible periodic Nielsen paths of $g'$ are fixed, the same is true for each $g_n$. Since $g_N=g$ has at most one indivisible fixed Nielsen path, it has at most one indivisible periodic Nielsen path, which if it exists is fixed. This proves items~(\ref{ItemNoNielsen}) and~(\ref{ItemRho}).

Now we prove item~(\ref{ItemGeneration}). By inductively applying Stallings fold factorization to the map $g \from G \to G$ we obtain a sequence of maps
$$G_0 \xrightarrow{f_1} G_1 \xrightarrow{f_2} \cdots \xrightarrow{f_K} G_K \xrightarrow{h} G_0
$$
where each $f_i$ is an edge isometric fold that preserves marking, and $h$ is a homothety. For $k=1,\ldots,K-1$ define $g_k \from G_k \to G_k$ to be 
$$g_k = f_k \composed \cdots \composed f_1 \composed h \composed f_K \composed \cdots \composed f_{k+1}
$$
and define $g_K \from G_K \to G_K$ to be $g_K = f_K \composed\cdots\composed f_1 \composed h$. Starting with $\rho_0=\rho$, for each $k=1,\ldots,K$ inductively define a path $\rho_k$ in $G_k$ by $\rho_k = (f_k \composed \rho_{k-1})_\# = (f_k \composed\cdots\composed f_1\composed \rho)_\#$. 

We prove inductively that for each $k=1,\ldots,K$, the map $g_k$ is a stable train track representative of $\phi$ with unique indivisible periodic Nielsen path $\rho_k$, and $g_k$ is obtained from $g_{k-1}$ by folding $\rho_{k-1}$ --- in other words, the turn folded by $f_k$ is the illegal turn taken by $\rho_{k-1}$. Assuming this is true for $k-1$, the turn of $G_{k-1}$ folded by the map $h\composed f_K \composed \cdots \composed f_{k} \from G_{k-1} \to G_0$ is an illegal turn for $g_{k-1}$, but by stability $g_{k-1}$ has a unique illegal turn, namely the illegal turn taken by the Nielsen path $\rho_{k-1}$, and so $g_k$ is obtained from $g_{k-1}$ by folding $\rho_{k-1}$. This implies that $g_k$ is a stable train track representative of $\phi$. Moreover, $\rho_k = (f_k \composed \rho_{k-1})_\#$ is an indivisible Nielsen path for $g_k$, and by stability $\rho_k$ is unique.

Since the Nielsen path $\rho_k = \alpha_k * \bar\beta_k$ is obtained from $\rho_{k-1} = \alpha_{k-1} * \bar\beta_{k-1}$ by folding initial oriented segments of $\bar\alpha_{k-1}$ and $\bar\beta_{k-1}$, we can describe the situation in the following manner. Noting that $\Length(\alpha_0)=\Length(\beta_0) = \frac{1}{2} \Length(\rho_0) = \Length(G)$, there exists a partition
$$0 < R_K \le R_{K-1} \le \cdots \le R_1 < R_0=\Length(G)
$$
such that 
\begin{align*}
\alpha_k &= f_k \composed \alpha_{k-1} \restrict [0,R_k] \\
\beta_k   &= f_k \composed \beta_{k-1} \restrict [0,R_k]
\end{align*}
and the oriented segments $\overline{\alpha_{k-1} \restrict [R_k,R_{k-1}]}$ and $\overline{\beta_{k-1} \restrict [R_k,R_{k-1}]}$ are the segments of $\overline\alpha_{k-1}$ and $\overline\beta_{k-1}$ that are folded by $f_k$. In other words, $f_k$ is generated by the relation $\alpha_{k-1}(t) \sim \beta_{k-1}(t)$ for $t \in [R_k,R_{k-1}]$. It now follows by induction that the map $f_k \composed \cdots \composed f_1 \from G_0 \to G_k$ is generated by the relation $\alpha_0(t)=\beta_0(t)$ for $t \in [R_k,\Length(G)]$: what one needs for the induction step is that for any composition of surjective maps $A \xrightarrow{\alpha} B \xrightarrow{\beta} C$, if the relation $\sigma$ on $A$ generates~$\alpha$, if the relation $\tau$ on $B$ generates $\beta$, and if the relation $\tau'$ on $A$ is mapped onto the relation $\tau$ by the map $\alpha\cross\alpha \from A \cross A \to B \cross B$, then the relation $\sigma\union\tau'$ generates $\beta\composed\alpha$. Since the maps $g=g_0$ and $f_K \composed \cdots \composed f_1$ differ by the homeomorphism $h_K$, it follows that $g$ is generated by the relation that generates $f_K \composed\cdots\composed f_1$, namely the relation $\alpha(t)\sim\beta(t)$ for $t \in [R_K,\Length(G)]$. Since $\rho = (g \composed \rho)_\#$, since $\Length(\rho)=2\Length(G)$, and since $g$ stretches path length by $\lambda$, a short calculation shows that $R_K = \Length(G)/\lambda$.
\end{proof}

We can now state the results that characterize geometricity of $T_+$ and parageometricity of $\phi$. These characterizations are restatements of results from \OuterLimits\ and \BH.

\begin{proposition}
\label{PropGeometric}
Suppose that $\phi \in \Out(F_r)$ is fully irreducible with attracting tree $T_+$, and let $g \from G \to G$ be a train track representative of a positive power of $\phi$ satisfying case~(\ref{ItemNoNielsen}) or~(\ref{ItemNielsenUnique}) of Proposition~\ref{PropNielsenUnique}. 
\begin{enumerate}
\item 
\label{ItemTNongeometric} If $g$ satisfies case~(\ref{ItemNoNielsen}) of Proposition~\ref{PropNielsenUnique} --- if $g$ has no periodic Nielsen paths --- then the tree $T_+$ is nongeometric.
\item 
\label{ItemTGeometric} If $g$ satisfies case~(\ref{ItemNielsenUnique}) of Proposition~\ref{PropNielsenUnique} --- if $g$ is Nielsen unique --- then the tree $T_+$ is geometric. 
\end{enumerate}
When $g$ is Nielsen unique, letting $\rho \from [0,2\Length(G)] \to G$ be the Nielsen path, we have:
\begin{enumeratecontinue}
\item 
\label{ItemPhiGeometric}
$\phi$ is geometric if and only if $\rho$ is a closed path. In this case $\rho$ traverses every edge of $G$ exactly twice.
\item 
\label{ItemPhiParageometric}
$\phi$ is parageometric if and only if $\rho$ is not a closed path. In this case $\rho$ traverses some edge of $G$ exactly once, and $\rho$ traverses some other edge of $G$ at least thrice.
\end{enumeratecontinue}
\end{proposition}

The proofs of items~(\ref{ItemTNongeometric}) and~(\ref{ItemTGeometric}) will be carried out in the remainder of Section~\ref{SectionCharacterizing}. First, in Section~\ref{SectionDirectLimits}, we give a direct limit construction which produces the tree $T_+$ out of the train track map $g$; this leads quickly to a proof of item~(\ref{ItemTNongeometric}) in Section~\ref{SectionNoNielsen}. Then, in Section~\ref{SectionWedgeModel}, we describe the wedge model which applies to the Nielsen unique case; this leads to the proof of  item~(\ref{ItemTGeometric}) in Section~\ref{SectionGeometricProof}. 

We start with:

\begin{proof}[Proof of items (\ref{ItemPhiGeometric}) and (\ref{ItemPhiParageometric}), assuming (\ref{ItemTNongeometric}) and (\ref{ItemTGeometric})] Let $g \from G \to G$ be a Nielsen unique train track representative of a positive power of $\phi$. Item~(\ref{ItemPhiGeometric}) was proved in \BH. Combined with item~(\ref{ItemTGeometric}), it immediately follows that $\phi$ is parageometric if and only if $\rho$ is not closed. 

Suppose that $\phi$ is parageometric. If $\rho$ traverses every edge of $G$ at least twice then using $\Length(\rho)=2 \Length(G)$ it follows that $\rho$ traverses every edge exactly twice, but that implies that $\rho$ is closed, a contradiction. It follows that $\rho$ traverses some edge $E$ at most once. On the other hand, by Proposition~\ref{PropNielsenUnique} item~(\ref{ItemNielsenUnique}) the path $\rho$ traverses each edge of $G$ at least once, and so $\rho$ traverses $E$ exactly once. Using again that $\Length(\rho)=2 \Length(G)$ it follows that $\rho$ traverses some other edge at least thrice. This proves item~(\ref{ItemPhiParageometric}).
\end{proof}

\subsection{Direct limits} 
\label{SectionDirectLimits}

Let $\phi \in \Out(F_r)$ be fully irreducible with attracting tree $T_+$, and let $g \from G \to G$ be any train track representative of $\phi$. We describe here a method for constructing $T_+$ from $g$. From a more well known point of view, $T_+$ is the Gromov--Hausdorff limit of the universal covering trees $T_0,T_1,T_2,\ldots$ of the sequence of marked graphs $G=G_0,G_1,G_2,\ldots$ given inductively by $[G_i] = [G_{i-1}]\phi$. Here we shall recast this point of view, making use instead of direct limits. In this section we give a preliminary description of the direct limit of $T_i$ in the category of semimetric spaces and distance nonincreasing maps; this direct limit is denoted $T^\#$. Then, in Section~\ref{SectionNoNielsen} we study the case that $g$ satisfies item~(\ref{ItemNoNielsen}) of Proposition~\ref{PropNielsenUnique}, and in Sections~\ref{SectionWedgeModel} and~\ref{SectionGeometricProof} we study the case that $g$ satisfies item~(\ref{ItemNielsenUnique}); in both cases we will verify that the metric space associated to the semimetric space $T^\#$ is the direct limit of $T_i$ in the category of metric spaces and distance nonincreasing maps, and we will use this information to identify the metric space direct limit with $T_+$.

Consider the sequence in $\X_r$ defined by $[G_0]=[G]$, $[G_i]=[G_{i-1}]\phi$. The underlying unmarked and unmetrized graphs $G=G_0,G_1,G_2,\ldots$ are all identical, and the map $g \from G \to G$ is rewritten as $g_i \from G_i \to G_{i+1}$. We then define a marking $R_r \to G_i$, given inductively by postcomposing the marking $R_r \to G_{i-1}$ with the map $g_{i-1}$. It follows that each map $g_i \from G_i \to G_{i+1}$ respects markings (up to homotopy). We also define a metric on $G_i$, given inductively as the unique metric such that $g_{i-1}$ maps each edge of $G_{i-1}$ locally isometrically to $G_{i}$. It follows that $\Length(G_{i+1}) = \Length(G_i) / \lambda(\phi)$. There is a homothety $h_i \from G\to G_i$ that compresses Lebesgue measure by a factor of $\lambda(\phi)^{-i}$ such that the following composition equals $g^i$:
\begin{equation}
\label{GraphSequence}
G=G_0 \xrightarrow{g_0} G_1\xrightarrow{g_1}\cdots \xrightarrow{g_{i-1}} G_i\xrightarrow{h_i^\inv} G
\end{equation}
Note that since the maps $g_0,\ldots,g_{i-1}$ all respect marking, the map $g_{i-1} \composed\cdots\composed g_0 = h_i \composed g^i \from G \to G_i$ also respects marking, and so the homothety $h_i^\inv \from G_i \to G$ and the map $g^i \from G \to G$ change marking by exactly the same outer automorphism, namely~$\phi^i$.

Let $T_i$ be the universal cover of $G_i$, so we may regard $T_i$ as an $F_r$-tree. Choose lifts of the maps $g_i$ to obtain a sequence of surjective equivariant maps
\begin{equation}
\label{DirectSequence}
T_0 \xrightarrow{\tilde g_0} T_1 \xrightarrow{\tilde g_1} T_2 \xrightarrow{\tilde g_2} \cdots
\end{equation}
such that $\tilde g_i$ maps each edge of $T_i$ isometrically onto an arc of $T_{i+1}$. Each map $\tilde g_i  \from T_i \to T_{i+1}$ is therefore distance nonincreasing. Since $[T_i] = [T_{i-1}]\phi$ as points of $\X_r$, and since $T_+$ is the attracting point in $\overline\X_r$ of each forward orbit of $\phi$ on $\X_r$ \BookZero, it follows that $\lim_{i \to \infinity} T_i = T_+$ in $\overline\X_r$.

Let $T^\#$ denote the \emph{direct limit} of the sequence~(\ref{DirectSequence}). Set theoretically, this is the set of equivalence classes of the disjoint union of $T_0,T_1,T_2,\ldots$, where $x_i \in T_i$ is equivalent to $x_j \in T_j$ if there exists $k \ge i,j$ such that $x_i,x_j$ have the same image in $T_k$, that is, $\tilde g_{k-1} \composed \cdots \composed \tilde g_i(x_i) = \tilde g_{k-1} \composed\cdots\composed\tilde g_j(x_j)$. Let $[x_i] \in T^\#$ denote the equivalence class of $x_i \in T_i$. By surjectivity of the maps $\tilde g_i$, for each $i$ each equivalence class has the form $[x_i]$ for (at least) some $x_i\in T_i$. Define a semimetric on $T^\#$ by the formula $d_{T^\#}([x_i],[y_i]) =\lim_{i \to \infinity} d_{T_i}(x_i,y_i)$ where by induction we define $x_i = \tilde g_{i-1}(x_{i-1})$ and similarly for $y_i$; the formula is clearly well-defined independent of the choice of a representative in $T_0$ of a given point in $T^\#$, and the limit exists because $d_{T_i}(x_i,y_i)$ is a nonincreasing sequence of positive numbers. The actions of $F_r$ on the trees $T_i$ induce an isometric action of $F_r$ on $T^\#$. The map $x_i \to [x_i]$ is an $F_r$-equivariant surjective function $q^\#_i \from T_i \to T^\#$. This argument shows that $T^\#$ is in fact the direct limit of the sequence $T_0 \xrightarrow{\tilde g_0} T_1 \xrightarrow{\tilde g_1} \cdots$, in the category of semimetric spaces with an isometric $F_r$-action and distance nonincreasing maps which are $F_r$-equivariant.

In order to understand $T^\#$ more precisely we need a result from \BookOne.  Fix $i\ge 0$ and an immersed arc or circle $\gamma_i$ in $G_i$, and inductively define $\gamma_{j+1} = (g_{j} \composed \gamma_{j})_\#$ for $j\ge i$. 
As long as $\gamma_j$ is a nondegenerate path, the number of illegal turns in $\gamma_j$ is nondecreasing as a function of $j$. It follows that \emph{either} $\gamma_j$ is eventually degenerate \emph{or} the number of illegal turns in $\gamma_j$ eventually stabilizes.

\begin{lemma}[Lemma 4.2.6 of \BookOne] 
\label{LemmaLegalConcatenation}
With the notation as above, assume that $\gamma_j$ is not eventually degenerate, and choose $J$ so that for $j \ge J$ the number of illegal turns in $\gamma_j$ is constant. Then for each $j \ge J$ the immersion $\gamma_j$ is a legal concatenation of legal paths and Nielsen paths. \qed 
\end{lemma}

By a ``legal concatenation'' we mean that the turn at each concatenation point is a legal turn. 

\subsection{No Nielsen path, strong convergence, and the proof of~(\ref{ItemTNongeometric})} 
\label{SectionNoNielsen}

We are now in a position to prove item~(\ref{ItemTNongeometric}) of Proposition~\ref{PropGeometric}. Let $\phi \in \Out(F_r)$ be fully irreducible and suppose that $g \from G \to G$ is a train track representative of a positive power of $\phi$ such that $g$ has no periodic Nielsen path. 

Using Lemma~\ref{LemmaLegalConcatenation} together with the nonexistence of Nielsen paths it follows that, given $i$ and a sequence of immersed arcs or circles $\gamma_j$ in $G_j$ defined for all $j \ge i$ such that $\gamma_{j+1}=(g_{j} \composed \gamma_{j})_\#$, there exists $J$ such that \emph{either} $\gamma_j$ is legal for all~$j\ge J$, \emph{or} $\gamma_j$ is degenerate for all $j \ge J$. This has several consequences.

First we show that the semimetric space $T^\#$ described in Section~\ref{SectionDirectLimits} is actually a metric space. Consider two distinct points of $T^\#$, represented by a pair of points in $T_i$ that are connected by a geodesic $\gamma_i$. Inductively define  $\gamma_{j+1} = (g_{j} \composed \gamma_{j})_\#$ for $j\ge i$. Since the endpoints of $\gamma_i$ map to distinct points in $T^\#$, the path $\gamma_j$ is not eventually degenerate, and so it is eventually legal. It follows that $\Length(\gamma_j)$ is a positive constant for $j \ge J$, and this constant equals the semimetric distance in $T^\#$ between the two given points of $T^\#$. In other words, the distance between an arbitrary pair of distinct points in $T^\#$ is positive, so $T^\#$ is a metric space.

Next, the metric space $T^\#$ is an $\reals$-tree, because the defining conditions for an $\reals$-tree metric are closed conditions \cite{CullerMorgan:Rtrees}. Moreover, $T^\#$ is a minimal $\reals$-tree. To see why, it suffices to prove that each point of $T^\#$ lies on a bi-infinite geodesic. Consider a point of $T^\#$ represented by $x_i \in T_i$.
The image of $x_i$ downstairs in $G_i$ lies on some leaf of the expanding lamination of $\phi$, because the expanding lamination realized in $G_i$ projects surjectively to $G_i$. Lifting this leaf we obtain a bi-infinite legal geodesic in $T_i$ containing $x_i$, and the image of this geodesic in $T^\#$ is a bi-infinite geodesic containing $[x_i]$. Since $[x_i]$ is arbitrary in $T^\#$, this proves minimality of $T^\#$.

Next, by applying Lemma~\ref{LemmaLegalConcatenation} to loops in $G_i$ it follows that the sequence of translation distance functions for the $F_r$-trees $T_i$ converges to the translation distance function on $T^\#$, from which it immediately follows that $T^\# = T_+$. 

Finally, applying Lemma~\ref{LemmaLegalConcatenation} again to the geodesic $\gamma_i$ between any two points in $T_i$, if these two points map to distinct points in $T^\#$ then for sufficiently large $j$ the path $\gamma_j$ is legal and so embeds isometrically in $T^\#=T_+$ under $q^\#_i$. But this is precisely the definition of \emph{strong convergence} of the sequence of $F_r$-trees $T_i$ to the $F_r$-tree $T_+$. By the main result of \cite{LevittPaulin:geometric} it follows that the tree $T_+$ is not geometric.

\subsection{Definition of the wedge model and its stable foliation $\fol^s$} 
\label{SectionWedgeModel}

Suppose now that that $g \from G \to G$ is a Nielsen unique train track representative of (a positive power of) a fully irreducible $\phi \in \Out(F_r)$. Let $\rho \from [0,2L] \to G$ be the unique Nielsen path of $g$, where $L=\Length(G)$, and let $\rho = \alpha * \bar\beta$ be the legal decomposition of $\rho$.

The \emph{wedge model} is an extension of the train track map $g \from G \to G$ to a homotopy equivalence $k\from K \to K$, where $K$ is the \nb{2}complex obtained from $G$ by attaching a disc $W$ to $G$, identifying an arc on the boundary of $W$ with the path $\rho$. The description of $k$ requires imposing on $W$ the structure of a ``wedge''. From this structure we will also obtain a measured foliation on $K$ denoted $\fol^s$, called the \emph{stable foliation} of $k$. The action of $k$ preserves leaves of $\fol^s$ and multiplies the transverse measure by $\lambda(\phi)$. In Section~\ref{SectionGeometricProof} we will exhibit geometricity of $T_+$ by proving that $T_+$ is dual to $\fol^s$.

Choose locally isometric parameterizations $\alpha,\beta\from[0,L]\to G$. The wedge $W = \triangle ABC$ is the triangle in $\reals^2$ with vertices $A=(0,+1)$, $B=(0,-1)$, $C=(L,0)$. The attaching maps are $(x,y) \mapsto\alpha(x)$ for $(x,y) \in\overline{AC}$, and $(x,y)\mapsto \beta(x)$ for $(x,y) \in \overline{BC}$. The \nb{2}dimensional cell complex $K$ is obtained by attaching $W$ to $G$ in this manner. The \nb{1}skeleton of $K$ is equal to $G$ union the base $\overline{AB}$ of $W$, with identifications $A\sim \alpha(0)$ and $B\sim\beta(0)$. The \emph{dihedral valence} in $K$ of a \nb{1}cell $E$ of $K$ is the total number of times that $E$ is traversed by the attaching maps of the \nb{2}cells of~$K$; by definition $E$ is a free edge if its dihedral valence equals~1. Note that $\overline{AB}$ is a free edge of $K$. Also, any edge $E$ of $G$ has dihedral valence equal to the number of times that $\rho$ traverses $G$. By collapsing $W$ from the free edge $\overline{AB}$ we obtain a deformation retraction of $K$ onto $G$, and so we may regard the \nb{2}complex $K$ as being marked by the homotopy equivalence $R_r \to G \inject K$, where $R_r \to G$ is the given marking of the marked graph $G$. We may therefore identify $\pi_1(K) \approx \pi_1(G)\approx \pi_1(R_r) \approx F_r$. 

The measured foliation $\fol^s$ is induced by the vertical measured foliation on $W$ equipped with the transverse measure $\abs{dx}$. To check compatibility along the \nb{1}skeleton, observe that for each edge $E$ of $G$, among the segments of $\overline{AC} \union\overline{BC}$ that map onto $E$, the measures obtained on $E$ by pushing forward $\abs{dx}$ via the attaching map all agree with the usual Lebesgue measure on $E$.

Now we define the extension $k \from K \to K$ of $g \from G \to G$. Subdivide $W = P \union W'$ where $P$ is the subtrapezoid of $W$ with one base $\overline{AB} = W \intersect \{x = 0\}$ and with parallel base $W\intersect\{x = L / \lambda(\phi) \}$, and $W'$ is the subwedge $W'= \overline{W-P}$. We call $W'$ the \emph{collapsed subwedge} of $W$ for reasons about to become apparent. Define $k$ to take $P$ onto $W$, stretching the $x$-coordinate by $\lambda(\phi)$, and for $0 \le x_0\le L / \lambda(\phi)$ contracting the $y$-coordinate of vertical segment $P\intersect\{x=x_0\}$ by a factor of $c(x_0)$, where $c(x)$ is the unique affine function satisfying $c(0)=1$, $c(L/\lambda(\phi)) = 0$. In~$W'$, for $L / \lambda(\phi) \le x_0 \le L$ the vertical segment $W' \intersect \{x=x_0\}$ is mapped by $k$ to the point $g(\alpha(x_0)) = g(\beta(x_0))$; the latter equation follows from the definition of an indivisible Nielsen path. This completes the definition of the wedge model $k\from K \to K$. 

We note a few facts about $k$ which will be important in what follows, and which are simple consequences of item~(\ref{ItemNielsenUnique}) of Proposition~\ref{PropNielsenUnique}:
\begin{itemize}
\item $k(x)=g(x)$ for each $x \in G$, so $k$ is continuous.
\item $k$ is a homotopy equivalence of $K$, and the induced outer automorphism on $\pi_1(K) \approx
F_r$ is $\phi$.
\item For each $x \in K$, the set $k^\inv(x)$ is described as follows:
\begin{itemize}
\item If $x \in K-G$ then $k^\inv(x)$ is a point in $K-G$.
\item If $x \in G$ then $k^\inv(x)$ is either a point in $G$ or a finite, connected graph whose edges are
a (uniformly) finite union of vertical segments of the subwedge $W'$.
\item Every vertical segment of $W$ is eventually collapsed by some power of $k$, \emph{except} for the
base $\overline{AB}$; this follows because in the Nielsen path $\rho = \alpha * \bar \beta$, with
$\alpha,\beta \from [0,L] \to G$, for each $t \in (0,L]$ there exists $i \ge 1$ such that
$g^i(\alpha(t))=g^i(\beta(t))$. 
\end{itemize}
\end{itemize} 
To understand the description of $k^\inv(x)$, note that the map $k$ is defined by collapsing to a point each of the vertical segments of the subwedge $W'$. The relation described in item~(\ref{ItemGeneration}) of Proposition~\ref{PropNielsenUnique}, which generates the map $g$, is exactly the same as the relation of endpoint pairs of vertical segments of the collapsed subwedge $W'$. It follows that for each $x \in G$,  the set $k^\inv(x)$ is a finite, connected union of vertical segments of $W'$: the vertices of this graph are the points of $g^\inv(x)$, a uniformly finite set, and each point of $G$ is an endpoint of a uniformly finite number of segments of $W'$, so $k^\inv(x)$ is a uniformly finite connected graph. The graph $k^\inv(x)$ can therefore be described by picking some vertical segment $\ell_0$ of $W'$, then inductively defining $\ell_i$ be the union of $\ell_{i-1}$ with all vertical segments of $W'$ that touch $\ell_{i-1}$, and then taking the union of the $\ell_i$ to obtain $k^\inv(x)$.

Let $K_0=K$ and let $K_i$ be the marked \nb{2}complex similarly obtained from the marked graph $G_i$ by attaching a wedge along the Nielsen path. The homotopy equivalence $k \from K \to K$ induces a marking preserving homotopy equivalence $k_i \from K_i \to K_{i+1}$ that agrees with $g_i \from G_i \to
G_{i+1}$. Note that the homothety $h_i \from G \to G_i$ defined earlier extends to a homeomorphism also denoted $h_i \from K \to K_i$ such that the composition $K=K_0 \xrightarrow{k_0} K_1 \xrightarrow{k_1} \cdots \xrightarrow{k_{i-1}} K_i \xrightarrow{h_i^\inv} K$ equals $k^i$. Letting $\fol^s_0=\fol^s$, there is a measure foliation $\fol^s_i$ defined inductively on $K_i$ as the pushforward of $\fol^s_{i-1}$ via the map $k_{i-1}$. Note that under the homemorphism $h_i \from K\to K_i$, $\fol^s_i$ is the pushforward of $\fol^s$ with transverse measure multiplied by $\lambda(\phi)^{-i}$. Lifting to universal covers, we obtain a \nb{2}complex $\wt K_i$ containing the tree $T_i$, and an action of $F_r$ on $\wt K_i$ extending the action on $T_i$. The map $\wt g_i \from T_i\to T_{i+1}$ extends to a map $\wt k_i\from \wt K_i \to \wt K_{i+1}$ that is a lift of $k_i\from K_i \to K_{i+1}$. There is an $F_r$-equivariant measured foliation $\wt\fol^s_i$ on $\wt K_i$ that is the lift of $\fol^s_i$ as well as the pushforward of $\wt\fol^s_{i-1}$ via $\tilde k_{i-1}$. 

\subsection{Proof of~(\ref{ItemTGeometric}): duality of $\fol^s$ and $T_+$} 
\label{SectionGeometricProof}
Given a fully irreducible $\phi \in \Out(F_n)$ with attracting tree $T_+$, in this section we prove item~(\ref{ItemTGeometric}) of Proposition~\ref{PropGeometric}: that under the assumption of Nielsen uniqueness, $T_+$ is a geometric $F_r$-tree.

Recall the notation: $g \from G \to G$ is a Nielsen unique train track representative of $\phi^i$ for some $i>0$, with Nielsen path $\rho$. From the constructions of Section~\ref{SectionWedgeModel}, let $k \from K \to K$ be a wedge model extension of $g$, with wedge $W$ attached to $G$ along $\rho$ to obtain $K$. Let $\fol^s$ be the stable measured foliation of $k$ on the \nb{2}complex $K$, induced by the vertical measured foliation on $W$.

Recall also the sequence~(\ref{GraphSequence}) of marked graphs and marked homotopy equivalences $G=G_0 \xrightarrow{g_0} G_1 \to \cdots$, the universal cover sequence~(\ref{DirectSequence}) of $F_r$-trees and $F_r$-equivariant maps $T_0 \xrightarrow{\tilde g_0} T_1 \to\cdots$, and $T^\#$ the direct limit of this sequence  in the category of $F_r$-semimetric spaces and $F_r$-equivariant distance nonincreasing maps.

Let $\L$ denote the leaf space of the measured foliation $\wt\fol^s$, with $F_r$-equivariant semimetric $d_\L$ defined as follows. Recall the semimetric $d_{\wt K}(x,y)$ on $\wt K$, defined by integrating the transverse measure of $\wt\fol^s$ along paths connecting $x$ to $y$ and taking the infimum. When $x,y$ are in the same leaf of $\fol^s$, clearly $d_{\wt K}(x,y)=0$ by integrating along a leaf segment connecting $x$ to $y$. The semimetric $d_{\wt K}$ therefore induces a well-defined semimetric $d_{\L}$ on $\L$. Also, the $F_r$-map from $\wt K$ to the dual $\reals$-tree of $\fol^s$ factors as the composition of the natural quotient $F_r$-map $\wt K \to \L$ and a surjective $F_r$-map from $\L$ to the dual $\reals$-tree. 

We shall prove geometricity of $T_+$ by proving that $T_+$ is equivariantly isometric to the dual $F_r$-tree of $\fol^s$ --- the metric space associated to the semimetric $d_{\wt K}$ on $\wt K$. We proceed indirectly, working with the semimetric space $T^\#$ instead of the tree $T_+$, and with the leaf space $\L$ of $\wt\fol^s$ instead of the dual $R_r$-tree of $\fol^s$.

Define an $F_r$-equivariant surjective map $\alpha \from T^\# \to \L$ as follows: for $x \in T_0$ representing $[x] \in T^\#$, $\alpha[x]$ is the leaf of $\wt\fol^s$ passing through $x$. To see that $\alpha$ is well defined, note that for $x,y \in T_0$ we have $[x]=[y]$ if and only if $x,y$ have the same image in some $T_i$, which occurs if and only if there exists a sequence of vertical segments of the wedge $W$ connecting $x$ to $y$ each of which are eventually collapsed by $\tilde k$, but this implies that $x,y$ are in the same leaf of $\wt\fol^s$. 

We claim that:
\begin{itemize}
\item[(A)] The map $\alpha \from T^\# \to \L$ is distance preserving. 
\end{itemize}
To prove this, consider $x_0,x'_0 \in T_0$ contained in respective leaves $\ell,\ell' \in \L$ of the measured foliation $\wt\fol^s$, and let $\xi = [x_0]$, $\xi'=[x'_0] \in T^\#$, so $\alpha(\xi)=\ell$ and $\alpha(\xi')=\ell'$. Define inductively $x_i = \tilde g_{i-1}(x_{i-1})$ and similarly for $x'_i$. Let $\gamma_i$ be the geodesic in $T_i$ between $x_i$ and $x'_i$, so the sequence $\Length_{T_i}(\gamma_i)$ is nonincreasing and has limit $d_{T^\#}(\xi,\xi')$. Choosing $\epsilon>0$, for sufficiently large $I$ we have $\Length_{T_I}(\gamma_I) \le d_{T^\#}(\xi,\xi') + \epsilon$. Starting with $\rho_I = \gamma_I$, inductively define a path $\rho_i$ in $K_i$ for $i=I-1,\ldots,0$ as follows: at any point where $\rho_{i+1}$ does not pull back continuously to $\wt K_{i}$, one can interpolate a leaf segment of $\wt\fol^s_{i}$, producing a continuous path $\rho_{i}$ in $\wt K_{i}$ connecting $x_i$ to $x'_i$, such that $\int_{\rho_{i}} \wt\fol^s_{i} = \int_{\rho_{i+1}} \wt\fol^s_{i+1}$. We therefore have 
$$\int_{\rho_0} \wt\fol^s_0 = \int_{\rho_I} \wt\fol^s_I = \Length_{T_I}(\gamma_I)
< d_{T^\#}(\xi,\xi')+\epsilon
$$
Letting $\epsilon\to 0$, it follows that $d_\L(\alpha(\xi),\alpha(\xi')) \le
d_{T^\#}(\xi,\xi')$. For the opposite inequality, choosing $\epsilon>0$ let $\rho_0$ be a continuous path
in $\wt K$ from $x_0$ to $x'_0$ so that $\int_{\rho_0} \wt\fol^s \le d_T(\ell,\ell') + \epsilon/2$.
Without increasing the integral along $\rho_0$ we may rewrite it as a concatenation of immersed paths in
$T_0$ and vertical segments of wedges. Inductively define the path $\rho_i$ in $\wt K_i$ as $\wt k_{i-1}
\composed\rho_{i-1}$, which has the effect of collapsing certain vertical wedge segments of $\rho_{i-1}$,
and so $\Length(\rho_i)=\Length(\rho_{i-1})$. By induction, $\Length(\rho_i)=\Length(\rho_0)$. For
sufficiently large $i$, say $i \ge I$, all vertical wedge segments of $\rho_0$ have been collapsed in
$\rho_i$ except for those which are lifts of the base $\overline{AB}$ of the wedge
$W$; let $\gamma_i$ be the path in $T_i$ obtained from $\rho_i$ be replacing each such vertical segment
with the associated Nielsen path. For $i \ge I$ the number of these Nielsen paths is constant, and
their length goes to zero as $i \to \infinity$, and so for sufficiently large $i$ we have
$$d_{T^\#}(\xi,\xi')\le\Length(\gamma_i) \le \Length(\rho_i) + \epsilon/2 < d_T(\ell,\ell') + \epsilon
$$
This proves the claim.

We claim next that:
\begin{itemize}
\item[(B)] For each $\xi,\xi' \in T^\#$, $\alpha(\xi) = \alpha(\xi')$ if and only if $d_{T^\#}(\xi,\xi')
= 0$.
\end{itemize}
If $d_{T^\#}(\xi,\xi') \ne 0$ then the previous claim shows that $d_\L(\alpha(\xi),\alpha(\xi')) \ne 0$
and so $\alpha(\xi) \ne \alpha(\xi')$. To prove the converse, suppose that $d_{T^\#}(\xi,\xi')=0$.
Choosing $x_0,x'_0 \in T_0$ so that $\xi=[x_0]$, $\xi'=[x'_0]$, we must prove that the points $x_0,x'_0$
are contained in the same leaf of $\wt\fol^s$. Define inductively $x_i=\tilde g_{i-1}(x_{i-1})$ and
similarly for $x'_i$. Let $\gamma_i$ be the geodesic between $x_i$ and $x'_i$ in $T_i$. Applying
Lemma~4.2.6 of \BookOne\ (Lemma~\ref{LemmaLegalConcatenation} above), for sufficiently large $i$, say $i
\ge I$, the path $\gamma_i$ is a legal concatenation of legal paths and Nielsen paths. As $i \ge I$
increases, the lengths of the legal paths in $\gamma_i$ stay the same while the number of Nielsen paths
is constant and their lengths go to zero. It follows that $d_{T^\#}(\xi,\xi') = 0$ only if $\gamma_i$ has
no legal paths for $i \ge I$, that is, $\gamma_I$ is a concatenation of Nielsen paths. We may therefore
connect $x_I$ to $x'_I$ by a path $\rho_I$ in $\wt K_I$ entirely contained in a leaf of $\wt\fol^s_I$.
Now proceeding inductively as in the earlier claim, for $i=I,I-1,\ldots,0$ we obtain a path $\rho_i$ in
$\wt K_i$ entirely contained in a leaf of $\wt\fol^s_i$ connecting $x_i$ to $x'_i$, and taking $i=0$ it
follows that
$\alpha(\xi)=\alpha(\xi')$. 

From the above two claims we can draw the conclusion that the semi-metric on $\L$ is a metric: if
$\ell,\ell'\in \L$ and $d_{\L}(\ell,\ell')=0$ then, choosing $\xi \in \alpha^\inv(\ell)$ and $\xi' \in
\alpha^\inv(\ell')$, Claim (A) implies that $d_{T^\#}(\xi,\xi')=0$, and Claim (B) implies that
$\ell=\alpha(\xi)=\alpha(\xi')=\ell'$. It immediately follows, follows from the definition of $\L$, that
$\L$ is equal to the associated metric space of the semimetric on $\wt K$, in other words, $\L$ is equal
to the dual tree of $\fol^s$. It also immediately follows that the map $\alpha$ identifies $\L$ with the
metric space associated to the semimetric on $T^\#$. Moreover, we know exactly which distinct pairs of
points in $T^\#$ have distance zero, namely, those pairs $\xi,\xi'$ represented by $x_0,x'_0 \in T_0$
contained in the same leaf $\ell$ of $\wt\fol^s_0$ so that $x_0$ and $x'_0$ are separated from each
other in $\ell$ by the lifts of $\overline{AB}$ in $\ell$. 

To complete the proof that $T_+$ is geometric, we must check that the translation distance function of the $F_r$-tree $\L$ is the limit of the translation distance functions of the $F_r$-trees $T_i$, for that will identify $\L$ with $T_+$. Consider $c \in \C$ represented by an immersed closed curve $\gamma_0$ in $G_0$. Inductively define $\gamma_i = (g_i\composed \gamma_{i-1})_\#$, so $\Length_{G_i}(\gamma_i)$ is the translation length of $c$ in $T_i$. By Lemma~\ref{LemmaLegalConcatenation}, for sufficiently large $i$, say $i\ge I$, $\gamma_i$ is a legal concatenation of legal paths and Nielsen paths, where the number of Nielsen paths is a constant independent of $i$ and their lengths go to zero, so $\lim\bigl(\Length_{G_i}(\gamma_i)\bigr)$ is equal to the total length of the portion of $\gamma_i$ which is not in one of the Nielsen paths, the latter number being independent of $i$. Let $A_i \subset T_i$ be an axis of (any representative of) $c$ acting on $T_i$. The map $A_i \to T^\# \to \L=T_+$ has the effect of folding each Nielsen path in
$A_i$ into a segment in $T^\#$, and mapping the rest of $A_i$ onto the axis $A_+$ of $c$ in $T_+$. It follows that a
fundamental domain for $A_+$ has length equal to the length of a fundamental domain for $A_i$ minus the Nielsen paths, which equals $\lim\bigl(\Length_{G_i}(\gamma_i)\bigr)$.

This completes the proof of Proposition~\ref{PropGeometric}.

\paragraph{Local topology of the wedge model.} When a fully irreducible outer automorphism $\phi \in \Out(F_r)$ has a geometric attracting tree $T_+$, there are some interesting connections between the behavior of $\phi$ and the topology of the wedge model $K$. For instance, when $\phi$ is geometric then $K$ is a surface with one boundary component; this is proved in \BH. 

When $\phi$ is parageometric, we shall need the following fact, an immediate consequence of Proposition~\ref{PropGeometric}~(\ref{ItemPhiParageometric}) and the construction of the wedge model. This fact will play a crucial role in the proof of Theorem~\ref{TheoremLambdaBigger}.

\begin{fact}
\label{FactParaWedgeModel}
If $\phi \in \Out(F_r)$ is parageometric, if $g \from G \to G$ is a Nielsen unique train track representative of a positive power of $\phi$, and if $k \from K \to K$ is the wedge model for $g$, then some edge of $G$ is a free edge of $K$, and some edge has dihedral valence~$\ge 3$.
\qed\end{fact}

\section{The stable foliation of the wedge model} 
\label{SectionStableFoliation}
In this section we fix a fully irreducible, parageometric outer automorphism $\phi \in \Out(F_r)$, a Nielsen unique train track representative $g \from G \to G$, and the wedge model $k \from K \to K$ of $g$, and we study the stable foliation $\fol^s$ of $k$, a $k$-invariant foliation of the \nb{2}complex $K$ whose leaves are compressed by the action of $k$. The main result, Proposition~\ref{PropHullLambdaMinus}, says that the set of bi-infinite lines in $\fol^s$, called the \emph{hull} of $\fol^s$, can be identified with the leaves of the expanding lamination of~$\phi^\inv$. The results of this section are closely related to results found in \OuterLimits. Indeed, we believe that Proposition~\ref{PropHullLambdaMinus} can be proved by the Rips machine methods of \OuterLimits, but we have developed a different proof.

\subsection{The stable foliation} For each $x \in K$ the leaf of $\fol^s$ through $x$ is called the \emph{stable leaf} of $x$, denoted $\fol^s_x$. We can build $\fol^s_x$ up inductively: let $\ell_0=x$, define $\ell_i$ inductively as the union of $\ell_{i-1}$ with all vertical segments of $W$ that intersect $\ell_{i-1}$, and then $\fol^s_x =\union_i \ell_i$. Note that $\fol^s_x$ is a locally finite, connected \nb{1}complex with vertex set $\fol^s_x \intersect G$ and whose edges are vertical segments of the wedge $W$. An edge path in the leaf $\fol^s_x$ is therefore a concatenation of a consecutive sequence of vertical segments of $W$, and the number of such segments is called the \emph{leafwise length} of the edge path. Any $p,q \in \fol^s_x\intersect G$ are connected by an edge path in $\fol^s_x$, and the smallest leafwise length of such a path is called the \emph{leafwise distance} between $p$ and~$q$. As we will see in Fact~\ref{FactTree} below, each leaf $\fol^s_x$ is a tree, so the minimal edge path $[p,q]$ between vertices $p,q \in \fol^s_x$ is unique, and the number of edges on this path is denoted $\Length_K [p,q]$, called the \emph{leafwise distance} between $p$ and $q$. Notice that we do \emph{not} measure $\Length_K$ using lengths of segments in the Euclidean triangle representation $W = \triangle ABC$. One should beware that the leafwise distance between $p$ and $q$ is not a continuous function, when regarded as a function on the set of ordered pairs $(p,q) \in G\cross G$ such that $p,q$ are contained a common leaf of $\fol^s_x$, because as $p,q$ vary in $G$ the edge path $[p,q]$ could vary in such a way that it
passes over the apex of $W$ where the leafwise distance jumps discontinuously. 

For each $x \in G$ the valence of the leaf $\fol^s_x$ at $x$ equals the number of times that
$\rho$ passes over $x$ except at an endpoint or the midpoint of the domain of $\rho$. If~$x$ is in the
interior of an edge $E$ of $K$, this number is just the dihedral valence of $E$ in $K$. Since there
are only finitely many vertices, it follows that there is a uniform upper bound for the valences of all
vertices in all leaves of $\fol^s_x$.

One point of confusion is the fact that, near the apex of the wedge $W$, the vertical segments
of $W$ get shorter and shorter, and it may seem possible that such segments could accumulate in a leaf
of $\fol^s_x$. This is not possible, however, because of local finiteness of the complex $K$: if $v$
is the vertex of $K$ to which the apex of $W$ is identified, and if $U$ is a small neighborhood of
$v$ in $K$, then for every vertical segment $\alpha$ of $W$ contained in $U$, and for every vertical
segment $\alpha'$ of $W$ such that $\alpha$ and $\alpha'$ share a common endpoint, the segment $\alpha'$
is not contained in $U$; see also the end of the proof of Lemma~\ref{KeyLemma}. In fact one sees that the
path topology on each leaf $\fol^s_x$, which has as basis the path components of $\fol^s_x \intersect
U$ over all open subsets $U \subset K$, is the same as the CW-topology on the simplicial complex
$\fol^s_x$. 

\renewcommand\ss{{ss}}

From the definition of $\fol^s$ and the map $k$ it is clear that $k$ preserves the foliation~$\fol^s$,
mapping each leaf $\fol^s_x$ onto the leaf $\fol^s_{k(x)}$, inducing a bijection of leaves. The next fact (almost) justifies the terminology ``stable foliation'' for $\fol^s$. Define the
\emph{strong stable set} of $x \in K$ to be
$$\fol^\ss_x = \{y \in K \suchthat \exists i \ge 0 \,\, \text{s.t.}\, k^i(x)=k^i(y)\}
$$

\begin{fact}[The stable foliation]
\label{FactStable}
Two points $x,y \in K$ are in the same leaf of $\fol^s$ if and only if there exists $i \ge 0$ such
that $k^i(x)=k^i(y)$ or $k^i(x),k^i(y) \in \overline{AB}$. It follows that for each $x \in K$ we have:
\begin{enumerate}
\item If $\overline{AB} \not\subset \fol^s_x$ then $\fol^\ss_x = \fol^s_x$.
\item If $\overline{AB} \subset \fol^s_x$ then the following hold. 
\begin{enumerate}
\item If $x \not\in
\interior(\overline{AB})$ then $\fol^\ss_x$ is the component of $\fol^s_x - \interior(\overline{AB})$
that contains~$x$. 
\item If $x \in \interior(\overline{AB})$ then $\fol^\ss_x = \{x\}$.
\end{enumerate}
\end{enumerate}
\end{fact}

\begin{proof}
This follows immediately from the observation that for any vertical segment $\gamma$ of the wedge
$W$, there is a power of $k$ collapsing $\gamma$ to a point if and only if $\gamma \ne \overline{AB}$. 
\end{proof}

\begin{fact}[Stable leaves are infinite] 
\label{FactInfiniteLeaves}
Each leaf of $\fol^s$ is an infinite \nb{1}complex.
\end{fact}

\begin{proof}
By Fact~\ref{FactStable}, for each $x \in G$ and each $n \ge 0$, the vertex set of the leaf of $\fol^s$
through $x$ contains the set $g^{-n}(g^n(x))$, whose cardinality goes to infinity as $n \to
\infinity$ since the transition matrix of $g$ is Perron--Frobenius. 
\end{proof}

\subsection{The structure of stable leaves} 

Consider the universal covering spaces $\wt G \subset \wt K$ of $G \subset K$. Let $\wt\fol^s$ be the foliation of $\wt K$ obtained by lifting $\fol^s$. 

Given a leaf $\ell$ of $\wt\fol^s$, we study the structure of $\ell$ by considering a locally embedded finite edge path $p$ in $\ell$. Since $p$ is an edge path, it starts and ends on $\wt G$.  Let the sequence of points of intersection of $p$ with $\wt G$ be denoted $x_0,x_1,\ldots,x_J$ where $J$ is the leafwise length of $p$. For example, one question we want to answer is whether $\ell$ is a tree, which is true if and only if $x_0 \ne x_J$ for all nontrivial edge paths $p$ in $\ell$; see Fact~\ref{FactTree}. Another question is whether $\ell$ is quasi-isometrically embedded in $\wt K$, which is true if and only if the map $i \mapsto x_i$ is a quasi-isometric embedding $\Z \to \wt K$; see Fact~\ref{FactQuasigeodesic}.

We establish further notation regarding $p$. For each $j=1,\ldots,J$, the segment of the path $p$ from
$x_{j-1}$ to $x_j$ is a vertical segment in a lift of the wedge, and so there exists a lift $\rho_j
\from [0,2 L] \to \wt G$ of $\rho$ or $\bar\rho$, and there exists $t_j \in [0,L)$, such that $x_{j-1} =
\rho_j(L-t_j)$ and $x_j = \rho_j(L+t_j)$. Let $V_j = \rho_j(L) \in \wt G$ be the point at which $\rho_j$
makes its unique illegal turn. The geodesic $\geodesic{x_{j-1}}{x_j}$ in $\wt G$ is a path of length
$2t_j$ obtained from $\rho_j$ by truncating the initial and final segments of length $L-t_j$, and so
$\geodesic{x_{j-1}}{x_j} = \alpha_j * \bar \beta_j$ concatenated at the illegal turn at $V_j$ with
$\alpha_j,\beta_j$ legal paths of length $t_j$. 

As a final remark, note in the discussion above that for each $j=1,\ldots,J-1$, since the subpath
of $p$ from $x_{j-1}$ to $x_{j+1}$ in $\ell$ is locally embedded, the vertical wedge segments from
$x_{j-1}$ to $x_j$ and from $x_j$ to $x_{j+1}$ are not inverses of each other, and so the paths $\alpha_j
* \bar\beta_j$ and $\alpha_{j+1} * \bar\beta_{j+1}$ are not inverses of each other in $\wt G$. 
Since $G$ has a unique illegal turn, each point of $\wt G$ at which an illegal turn
occurs has a unique illegal turn, from which it follows that $V_j \ne V_{j+1}$.

\begin{lemma} 
\label{KeyLemma} Let $J \ge 1$ and let $p$ be a locally embedded edge path of leafwise length~$J$ in a leaf of~$\wt\fol^s$. Using the notation above, for each $j = 0,\ldots,J$ there exists $i$ with $1 \le i \le J$ such that the geodesic $\geodesic{x_0}{x_J}$ makes an illegal turn at $V_i$, and the path $\geodesic{x_j}{V_i}$ is legal. Moreover, $\Length(\geodesic{x_j}{V_i}) \le L$.
\end{lemma}

\begin{proof} We will prove the final statement about length at the very last. The proof of the rest
of the lemma is by induction on $J$, with $J=1$ obvious. Supposing that the
lemma is true for a certain $J$, we wish to prove it for $J+1$. 

Choose $x_j$ with $0 \le j \le J+1$. We may assume that $j \le J$, because if $j=J+1$ then we
can just reverse the direction of~$p$, reducing to the case $j=0$. We therefore can apply the
induction hypothesis to the subpath of $p$ from $x_0$ to $x_J$, obtaining $k$ with $1 \le k \le J$ such
that $\geodesic{x_0}{x_J}$ has an illegal turn at $V_{k}$ and such that $\geodesic{x_j}{V_k}$ is legal.

Decompose $\geodesic{x_0}{x_J} = \mu * \nu$ at the turn $V_{k}$; no claims are made on the legality of
$\mu$ or $\nu$.  The path $\geodesic{x_0}{x_{J+1}}$ is what you get by tightening $\mu * \nu *
\alpha_{J+1} * \bar\beta_{J+1}$, which is done by cancelling a terminal segment of $\mu * \nu$ with an
initial segment of $\alpha_J * \bar\beta_J$.

If the cancellation does not remove all of $\nu$ then $V_{k}$ is still a turn in $\geodesic{x_0}{x_{J+1}}$ and we are done, by taking $V_i = V_k$. So we may assume that the terminal segment of $\mu*\nu$ which cancels is, at least, all of $\nu$. Since $\alpha_{J+1}$ is legal it does not cancel with any of $\mu$. 

If not all of $\alpha_{J+1}$ cancels with $\nu$ then $\geodesic{x_0}{x_{J+1}}$ takes the illegal
turn $V_{J+1}$. To show that the path $\geodesic{x_j}{V_{J+1}}$ is legal, this path is a concatenation
of the legal path $\geodesic{x_j}{V_{k}}$ with the path $\geodesic{V_{k}}{V_{J+1}}$ which is legal
because it is a subpath of $\alpha_{J+1}$, and by construction the turn at which these two paths are
concatenated is not the illegal turn at the point $V_{k}$, and so this turn is legal. We are therefore
done, taking $V_i = V_{J+1}$.

The remaining case is when $\nu = \bar\alpha_{J+1}$, which we show leads to a contradiction. In this case $V_{k}=V_{J+1}$ and so the illegal turns of $\geodesic{x_0}{x_J} = \mu * \nu$ and $\alpha_{J+1} * \bar \beta_{J+1}$ at the common concatenation point are the same, because $\wt G$ has at most one illegal turn at each vertex. Also, the Nielsen paths $\rho_{J+1}$, $\rho_{k}$ are equal up to orientation, because each is a lift of $\rho$ or $\bar\rho$ and they have the same illegal turn. This implies that $(J+1) - k \ge 2$, from the final remark just before Lemma~\ref{KeyLemma}; this is where we use the hypothesis that $p$ is locally embedded. If the orientations of the Nielsen paths $\rho_k$ and $\rho_{J+1}$ agree, then the points $x_k$ and $x_{J+1}$ lie on the same legal half of this Nielsen path, and so the geodesic $\overline{x_k x_{J+1}}$ is legal. But this geodesic falls under the induction hypothesis because $1 \le J+1-k \le J$, and so $\overline{x_k x_{J+1}}$ must have an illegal turn, a contradiction. If $\rho_k$, $\rho_{J+1}$ are oriented in the opposite direction, then the path $\overline{x_k x_J}$ is legal, which also contradicts the induction hypothesis, noting that $1 \le J-k \le J$. The case $\nu = \bar\alpha_{J+1}$ therefore cannot occur, completing the induction. 

\paragraph{Remark:} Note the following consequence of what we have proved so far: any geodesic in $\wt G$ whose endpoints lie on the same leaf of $\wt\fol^s$ has at least one illegal turn.

\bigskip

It remains to prove that the legal  path $\geodesic{x_j}{V_i}$ has length $\le L$. Switching orientation of $p$ if necessary, we may assume $j<i$. It suffices to show that the legal paths $\geodesic{x_j}{V_i}$ and $\geodesic{x_i}{V_i}$ have the same length, because 
$$\Length \geodesic{x_i}{V_i} \le \frac{1}{2} \Length(\rho_i) = \frac{1}{2} \Length(\rho) = L
$$
Note that $\geodesic{x_j}{x_i}$ is obtained by concatenating two legal paths in $\wt G$, namely $\geodesic{x_j}{V_i}$ and the reverse of $\geodesic{x_i}{V_i}$, and then cancelling. The geodesic $\overline{x_j x_i}$ therefore has at most one illegal turn, but the remark above shows it has exactly one illegal turn. A similar argument shows that for each $k \ge 1$, if $\geodesic{\tilde g^k(x_j)}{\tilde g^k(x_i)} = (\tilde g^k(\geodesic{x_j}{x_i}))_\#$ is nontrivial then it has exactly one illegal turn: it has at most one illegal turn because it is the straightened image of a path with exactly one illegal turn; but it has at least one illegal turn by the remark above. The proof breaks now into two cases, depending on whether $x_i$ and $x_j$ are identified by some power of $\tilde g$.

\subparagraph{Case 1: $\tilde g^k(x_j) = \tilde g^k(x_i)$ for some $k \ge 1$.} In this case the legal paths $\tilde g^k \left(\geodesic{x_j}{V_i}\right)$ and $\tilde g^k \left(\geodesic{x_i}{V_i}\right)$ are equal, so they have the same length, but $g^k$ stretches the path length of every legal path by the same factor $\lambda^k$, and so $\Length\geodesic{x_j}{V_i} = \Length\geodesic{x_i}{V_i}$. 

\subparagraph{Case 2: $\tilde g^k(x_j) \ne \tilde g^k(x_i)$ for all $k \ge 1$.} In this case the geodesic $\geodesic{\tilde g^k(x_j)}{\tilde g^k(x_i)}$ has exactly one illegal turn, for all $k \ge 1$. Applying Fact~\ref{FactStable}, it follows that for any sufficiently large $k$ the points $\tilde g^k(x_j)$ and $\tilde g^k(x_i)$ are endpoints of some lift of $\overline{AB}$, and hence the geodesic $\geodesic{\tilde g^k(x_j)}{\tilde g^k(x_i)}$ is the corresponding lift of $\rho$, which we denote~$\tilde\rho$. But $\geodesic{\tilde g^k(x_j)}{\tilde g^k(x_i)}$ is obtained by concatenating the legal path  $\tilde g^k \left(\geodesic{x_j}{V_i}\right)$ with the reverse of the legal path $\tilde g^k \left(\geodesic{x_i}{V_i}\right)$ and then canceling, and what is left over after the cancellation is $\tilde\rho$ which is a concatenation of two legal paths of equal length, one an initial segment of $\tilde g^k \left(\geodesic{x_j}{V_i}\right)$ and the other an initial segment of $\tilde g^k \left(\geodesic{x_i}{V_i}\right)$. This proves that the legal paths  $\tilde g^k \left(\geodesic{x_j}{V_i}\right)$ and $\tilde g^k \left(\geodesic{x_i}{V_i}\right)$ have the same length, and so $\Length\geodesic{x_j}{V_i} = \Length\geodesic{x_i}{V_i}$ as in Case~1.
\end{proof}

\subsection{Consequences of Lemma \ref{KeyLemma}}

The first consequence is:

\begin{fact}
\label{FactTree}
Each leaf of $\wt\fol^s$ is a tree.
\end{fact}

\begin{proof}
Lemma~\ref{KeyLemma} shows that for each locally embedded edge path in a leaf of $\wt\fol^s$, the geodesic in $\wt G$ with the same endpoints has at least one illegal turn, and so the geodesic is nondegenerate and the endpoints are distinct. This implies that the leaf has no loops and so is a tree.
\end{proof}

\begin{fact}
\label{FactTreeDownstairs}
Each leaf of $\fol^s$ is a tree.
\end{fact}

\begin{proof} Let $\ell$ be a leaf of $\fol^s$. By Fact~\ref{FactTree}, $\ell$ is $\pi_1$-injective in
$K$. If $\ell$ were not a tree then then $\ell$ would contain an embedded loop. By
Fact~\ref{FactStable} the image of this loop under a sufficiently high power of the homotopy
equivalence $k$ is either a point or the segment $\overline{AB}$, violating $\pi_1$-injectivity. 
\end{proof}

Consider now any $F_r$-equivariant proper geodesic metric on $\wt K$, and any leaf $\ell$ of $\wt \fol^s$,
with the simplicial structure on $\ell$ inducing a simplicial metric where each edge has length~1. Next
we essentially prove that the injection $\ell \to \wt K$ is a quasi-isometric embedding, with
quasi-isometry constants independent of $\ell$. What we actually need, and prove, is an interpretation
of this statement that takes place entirely in $\wt G$:

\begin{fact} 
\label{FactQuasigeodesic} For each leaf $\ell$ of $\wt\fol^s$, and each $x,y \in \ell \intersect \wt G$,
letting $x=x_0,\ldots,x_J=y$ be the points of $\ell \intersect \wt G$ from $x$ to $y$ in order, the map
$i\mapsto x_i$ is a quasigeodesic embedding of $\{0,1,\ldots,J\}$ into $\wt G$, with quasigeodesic
constants independent of $x,y$.
\end{fact}

\begin{proof} Let $\Length(\cdot)$ denote length in $\wt G$. Since $J$ is arbitrary, it suffices to prove
that $\Length(\geodesic{x_0}{x_J})$ is bounded above and below by an affine function of $J$. Since
$\Length(\geodesic{x_i}{x_{i+1}}) \le 2L$ we have $\Length(\geodesic{x_0}{x_J})\le 2L \cdot J$. 

Lemma~\ref{KeyLemma} gives a map $f$ from the set $\{x_0,\ldots,x_J\}$ to the set of
illegal turns of the geodesic $\geodesic{x_0}{x_J}$, such that $d(x_j,f(x_j)) \le L$. There is an
integer $\kappa \ge 1$ such that the $L$ neighborhood of each illegal turn $V_i$ intersects at most
$\kappa$ different edges. Lemma~\ref{KeyLemma} also implies that each edge contains at most one of the
points $\{x_0,\ldots,x_J\}$, and so the $L$ neighborhood of $V_i$ contains at most
$\kappa$ of the points $\{x_0,\ldots,x_J\}$. This shows that the map $f$ is at most $\kappa$-to-one,
and so its image has cardinality $\ge \frac{J+1}{\kappa}$. In other words, the geodesic
$\geodesic{x_0}{x_J}$ has at least $\frac{J+1}{\kappa}$ illegal turns. Letting $\eta>0$ be the minimum
length of an edge of $G$, it follows that 
\begin{align*}
\Length(\geodesic{x_0}{x_J}) &\ge (\frac{J+1}{\kappa}-1) \eta 
= \frac{\eta}{\kappa} J + \eta (\frac{1}{\kappa}-1)
\end{align*}
\end{proof}

\subsection{The hull of the stable foliation}

The \emph{hull} of a simplicial tree is the union of bi-infinite lines in the tree, and the hull of
$\fol^s$ is the union of the hulls of its leaves, a closed foliated subset of $K$ denoted
$\Hull(\fol^s)$. The hull of $\wt\fol^s$ is similarly defined, and $\Hull(\wt\fol^s)$ is equal to the
total lift of $\Hull(\fol^s)$. By Fact~\ref{FactQuasigeodesic}, each bi-infinite line in a leaf of
$\Hull(\wt\fol^s)$ intersects $\wt G$ in a quasigeodesic embedding $\Z \to \wt G$ and so has finite
Hausdorff distance from a bi-infinite line in $\wt G$. Our goal in
this section is to identify this collection of lines with a natural extension of the lamination
$\Lambda^u(\phi^\inv)$.

\begin{fact} 
\label{FactHullInfo}
$\Hull(\fol^s)$ is a nonempty, closed lamination, using the compact-open topology on maps
$\reals \to K$. For each leaf $\ell$ of $\Hull\fol^s$, $k(\ell)$ is also a leaf of $\Hull\fol^s$, and
this induces a bijection of the set of leaves of $\Hull\fol^s$. The leaf $\ell'$ that maps to $\ell$ is
the unique leaf contained in the set $k^\inv(\ell)$.
\end{fact}

\begin{proof} 
As we saw in Fact~\ref{FactInfiniteLeaves}, each leaf has infinitely many vertices. Also, each vertex
has finite valence --- in fact, the valence of each vertex of each leaf is uniformly bounded, by the
maximum cardinality of $\rho^\inv(x)$ for $x \in G$. It follows that each leaf contains an infinite ray,
and so there exists a sequence of longer and longer leaf segments $\alpha_i$ centered on a sequence of
points $x_i$. Any limit point $x$ of the sequence $x_i$ lies on a bi-infinite line in some leaf. This
shows simultaneously that $\Hull(\fol^s)$ is closed and that it is nonempty.

The map $k$ is finite-to-one on $G$, and it extends to $K$ by collapsing to a point each vertical
segment contained in the collapsed subwedge $W'$ of $W$. It follows that each point pre-image is a
connected subset of a leaf of $K$ consisting of a union of boundedly many edges in that leaf, and since
each leaf is a tree (by Fact~\ref{FactTreeDownstairs}), each point pre-image is a finite subtree. The
intersection of $\ell$ with a point pre-image is therefore a finite subarc. In other words, the effect
of $k$ on $\ell$ is to collapse to a point each of a pairwise disjoint collection of uniformly finite
subarcs of $\ell$, which implies that $k(\ell)$ is also a bi-infinite line in a leaf of $\fol^s$. Thus
$k$ induces a well-defined self-map of the set of leaves of $\Hull\fol^s$. This map is injective,
because for each leaf $\ell$ of $\Hull\fol^s$, the set $k^\inv(\ell)$ maps to $\ell$ by collapsing a
pairwise disjoint collection of uniformly bounded finite subtrees, and so $k^\inv(\ell)$ contains a
unique leaf $\ell'$ of $\Hull(\fol^s)$. Since $k(\ell')$ is contained in $\ell$ and equals a leaf of
$\Hull\fol^s$, it follows that $k(\ell')=\ell$, showing that the map is surjective.
\end{proof}

Now we set up notation for pushing leaves of $\Hull(\wt\fol^s)$ from $\wt K$ into the graph~$\wt G$.
Consider a leaf $\ell$ of $\Hull(\wt\fol^s)$. Picking a base point $x_0 \in \ell \intersect \wt G$ and an
orientation of $\ell$ determines an ordering of the set $\ell \intersect \wt G$, giving a
bi-infinite sequence that we denote $\bx(\ell) = (\ldots,x_{-1},x_0,x_1,\ldots)$. This sequence has the
following properties:
\begin{enumerate}
\item For each $i$, the geodesic $\geodesic{x_{i-1}}{x_i}$ in $\wt G$ is a subpath of a lift $\wt \rho =
\wt\alpha * \wt{\bar\beta}$ of the Nielsen path $\rho$, consisting of equal length terminal subpaths of
$\wt\alpha$ and $\wt{\bar\beta}$.
\item The sequence $\bx(\ell)$ is one-to-one.
\end{enumerate}
Conversely, any sequence $\bx = (\ldots,x_{-1},x_0,x_1,\ldots)$ satisfying (1) and (2) is equal to
$\bx(\ell)$ for some leaf $\ell$ of $\Hull(\wt\fol^s)$ and some base point and orientation of $\ell$.

Given a leaf $\ell$ of $\Hull(\wt\fol^s)$ and $\bx(\ell) = (\ldots,x_{-1},x_0,x_1,\ldots)$ as above,
by Fact~\ref{FactQuasigeodesic} the sequence $\bx(\ell)$ is a quasi-isometric embedding $\Z \to \wt G$
with uniform constants independent of $\ell$. It follows that $\bx(\ell)$ fellow travels a unique
bi-infinite line in $\wt G$ denoted $\bar\ell$. We thus obtain sublaminations of the geodesic
laminations of $\wt G$ and of~$G$:
\begin{align*}
\wt\Lambda_\Hull &= \{\bar \ell \suchthat \text{$\ell$ is a leaf of $\Hull(\wt\fol^s)$} \}
\subset \Lambda_{\wt G}
\\
\Lambda_\Hull &= \text{the projection of $\wt\Lambda_\Hull$ to $\Lambda_G$}
\end{align*}
For a sequence $\bx$ satisfying (1) and (2) above, the image of $\bx$ under $\tilde g$ also has the
structure of a sequence satisfying (1) and (2), and $\tilde g$ induces a quasi-isometry from $\bx$ to
$\wt g(\bx)$. It follows that $\tilde g$ induces a quasi-isometry $\bar\ell \to \bar\ell'$ well-defined
up to bounded distance. We may therefore set $\wt\phi(\bar\ell) = \bar\ell'$, which induces a map
$\Lambda_\Hull \to \Lambda_\Hull$ downstairs, that is, $\Lambda_\Hull$ is an invariant
sublamination of the action of $\phi$ on the geodesic lamination of $G$.

A leaf $\ell$ of the geodesic lamination of $F_r$ is \emph{birecurrent} if its realization in any
(equivalently, some) marked graph has the property that each finite subpath occurs infinitely often in
both ends of $\ell$. Equivalently, $\ell$ is contained in the set of limits points of each of its
ends, where the convergence takes place in the geodesic lamination. Compactness of $K$ implies that
$\Hull(\fol^s)$ is compact and so, as a consequence of the Hausdorff maximum principle, $\Hull(\fol^s)$
has a nonempty minimal sublamination. Each leaf of a minimal sublamination is birecurrent. We therefore
have shown:

\begin{fact} $\Lambda_\Hull$ contains a birecurrent leaf.
\qed\end{fact}

Here is the main result of this section. Following the proof we will strengthen it by sketching how to
identify $\Lambda_\Hull$ completely. Recall that $\Lambda^u(\phi)$ denotes the expanding lamination of
$\phi$ and $\Lambda^u(\phi^\inv)$ is the expanding lamination of $\phi^\inv$. In the next proposition we
consider the realizations of each of these laminations in the train track map $g \from G \to G$
representing $\phi$.

\begin{proposition}
\label{PropHullLambdaMinus}
$\Lambda_\Hull \supset \Lambda^u(\phi^\inv)$.
\end{proposition}

\begin{proof} First we show that there is a uniform bound to the length of any legal subpath of any
leaf of $\Lambda_\Hull$. To see why, fix a leaf $\bar\ell$ of $\Lambda_\Hull$, and let $\bx(\ell) =
(\ldots,x_{-1},x_0,x_1,\ldots)$. For each $i \ne j$, each point of the geodesic $\geodesic{x_i}{x_j}$
has distance $\le L$ from some $x_k$ with $i \le k \le j$, and by applying Lemma~\ref{KeyLemma} it
follows that $x_k$ has distance $\le L$ from some illegal turn of $\geodesic{x_i}{x_j}$, so the illegal
turns on $\geodesic{x_i}{x_j}$ are spaced no more than $4L$ apart. Any finite subsegment of
$\bar\ell$ of length $>4L$ is contained in some $\geodesic{x_i}{x_j}$ and so contains an illegal turn.

Recall from \BookOne\ that a leaf $\ell$ of the geodesic lamination is \emph{weakly attracted} to
$\Lambda^u(\phi)$ if any finite subpath of a leaf of $\Lambda^u(\phi)$ is a subpath of $g^n(\ell)$ for
some large~$n$. Since the subpath can be a legal segment of arbitrarily large length, it follows that
$g^n(\ell)$ contains legal paths of arbitrarily large length. But we have shown that if $\ell \in
\Lambda_\Hull$ then $g^n(\ell) \in \Lambda_\Hull$ and so there is an upper bound to the length of legal
subpaths. It follows that \emph{no} leaf of $\Lambda_\Hull$ is weakly attracted to $\Lambda^u(\phi)$.
Applying Theorem~6.0.1 of \BookOne, it follows that each birecurrent leaf of $\Lambda_\Hull$ is a leaf of
$\Lambda^u(\phi^\inv)$, so there exists at least one leaf of $\Lambda^u(\phi^\inv)$ that is contained in
$\Lambda_\Hull$. By minimality of $\Lambda^u(\phi^\inv)$ it follows that $\Lambda^u(\phi^\inv) \subset
\Lambda_\Hull$. \end{proof}

For completeness sake we give a description of $\Lambda_\Hull$, but with details of proof only
sketched lightly since we do not need this description for our present purposes. We'll assume that $\bar
g \from \bar G \to \bar G$ is a train track representative of $\phi^\inv$ which is either Nielsen unique
or has no Nielsen path at all. The \emph{extended expanding lamination} of $\phi^\inv$ is defined to be
the union of $\Lambda^u(\phi^\inv)$ with finitely many other leaves, as follows. At a periodic vertex $p$
of $\bar g$, any periodic direction $d$ determines a ray $r_d$ under iteration of $\bar g$, and any
two periodic directions $d,d'$ determine a leaf $\ell=r_d \union r_{d'}$ in the extended expanding
lamination. Also, if $\rho$ is the Nielsen path connecting two points $p$ and $q$, then for any periodic
directions $d,d'$ at $p,q$ respectively, distinct from the directions of $\rho$ at its endpoints,
$\ell = r_d\union\rho \union r_{d'}$ is a leaf in the extended expanding lamination. A study of the
stabilization algorithm of \BH\ shows that the extended expanding lamination is well-defined,
independent of the choice of $\bar g$. 

We claim that $\Lambda_\Hull$ is the extended expanding lamination of $\phi^\inv$. This is a consequence
of the proof of Theorem~6.0.1 of \BookOne\ applied to a train track representative $g
\from G\to G$ for $\phi$. This proof is laid out with Steps 1, 2, and 3. In the end of Step~2, one
considers a geodesic $\ell$, for example any leaf of $\Lambda_\Hull$, whose realization in $G$ is not
weakly attracted to $\Lambda^u(\phi)$, meaning that the sequence $(g^i \ell)_\#$ does not develop longer
and longer segments that are leaf segments of $\Lambda^u(\phi)$. What one shows in this situation is that
there exists an immersed loop $\gamma$ in $\bar G$ of uniformly bounded length so that the sequence 
$\gamma_i = (\bar g^i \gamma)_\#$ develops longer and longer segments that are contained in and exhaust
$\ell$, as realized in $\bar G$. Applying Lemma~\ref{LemmaLegalConcatenation}, for some $I$ the loop
$\gamma_I$ is a legal concatenation of legal paths and Nielsen paths, and it follows that for $i \ge I$
the loop $\gamma_i$ is a legal concatenation of a uniformly bounded number of segments each of which is
either a leaf segment of $\Lambda^u(\phi^\inv)$ or of one of the finitely many leaves added to make the
extended expanding lamination of $\phi^\inv$. One then sees that $\ell$ is exhausted by such segments,
implying that $\ell$ itself is a leaf of the extended expanding lamination.

\section{Proof of Theorem \ref{TheoremLambdaBigger}}
\label{SectionMainProof}

Recall the notation: $\phi \in \Out(F_r)$ is parageometric, with $g \from G \to G$ a Nielsen unique train track representative of some positive power of $\phi$, with wedge model $k \from K \to K$, with $\fol^s$ the stable foliation of $k$, and with $\Lambda_\Hull$ the sublamination of $\Lambda_G$ obtained by pushing leaves of $\Hull(\fol^s)$ into $G$. 

To prove that $\lambda(\phi^\inv) < \lambda(\phi)$ we will make use of an intermediate quantity, the asymptotic compression factor of $k$ acting on leaves of $\Hull(\fol^s)$. On the one hand we show that $\lambda(\phi)$ is a strict upper bound for this factor, by using a symbolic dynamics argument to interpret $\lambda(\phi)$ geometrically, and by exploiting Fact~\ref{FactParaWedgeModel} which says that some edge of $G$ is a free edge in $K$. On the other hand, we use the result from Section~\ref{SectionStableFoliation}, that $\Lambda_\Hull$ is a sublamination of $\Lambda^u(\phi^\inv)$, to show that the asymptotic compression factor is equal to $\lambda(\phi^\inv)$ on the nose. 

In this section we distinguish two measurements of length. First, we use $\Length_K(\alpha)$ to denote combinatorial length of an edge path $\alpha$ in a leaf of $\fol^s$ or $\wt\fol^s$, that is, the number of vertical wedge segments in $\alpha$. We also use $\Length_G$ to denote length of a path in the graph $G$ or $\wt G$.

By Fact~\ref{FactParaWedgeModel}, there is at least one edge of $G$ that is a free edge of the \nb{2}complex~$K$. Let $G_1$ be the subgraph of $G$ consisting of the nonfree edges, that is, the edges of dihedral valence $>1$ in~$K$. By Fact~\ref{FactParaWedgeModel} it follows that $G_1$ is a nonempty, proper subgraph. Let $I_1 = \{x \in G \suchthat k^n(x) \in G_1 \,\text{for all}\, n \ge 0\}$. Since $g$ has constant stretch factor $\lambda(\phi)$, since the transition matrix of $g$ has a positive power, and since $G_1$ is a nonempty, proper subgraph of $G$, it follows that $I_1$ is a Cantor set in~$G_1$.

\begin{fact} 
\label{FactHullValence}
If $\ell$ is a leaf of $\Lambda_\Hull$ then $\ell \intersect G \subset I_1$.
\end{fact}

\begin{proof} The hypothesis means that $\ell$ is a bi-infinite line in a leaf of $\fol^s$. This implies that $\ell \intersect G \subset G_1$, because for each point $x \in G$ that is not a vertex, the valence of $x$ in its $\fol^s$ leaf equals the dihedral valence of the edge containing $x$. Since $k^n(\ell)$ is also a bi-infinite line in a leaf of $\fol^s$, it follows that $k^n(\ell) \intersect G \subset G_1$. This being true for all $n$, it follows that $\ell\intersect G \subset I_1$.
\end{proof}

Our next fact gives the strict inequality that we shall need in proving that
$\lambda(\phi^\inv) < \lambda(\phi)$.

\begin{fact}
\label{FactLambdaPrime}
There is a number $\lambda' < \lambda(\phi)$ such that for each $x \in I_1$ we have
$$\limsup_{n \to \infinity} \frac{1}{n} \log \abs{g^{-n}(x) \intersect I_1} \le \log(\lambda')
$$
\end{fact}

\begin{proof}  Recall that $\TG$ is the transition graph of $G$. Let $\TG_1$ be the subgraph
obtained from $\TG$ by throwing away each vertex associated to an edge of $G$ of dihedral valence
$=1$, and any directed edge of $\TG$ incident to such a vertex. Let $M_1$ be the transition matrix of
$\TG_1$. If $x$ is contained in the interior of the $i^{th}$ edge of $G$, then $\abs{g^{-n}(x)
\intersect I_1}$ is the sum of the entries in the $i^{th}$ column of $M_1^n$.  In all cases
$\abs{g^{-n}(x) \intersect I_1}$ is bounded above by $\abs{M_1^n}$, the sum of all coefficients in
$M_1^n$. It therefore suffices to show that 
$$
\limsup_{n \to \infinity} \frac{1}{n} \log \abs{M_1^n} \le \log(\lambda') < 
\log(\lambda).
$$ 
Since we may regard $M_1$ as defined on the same set $\E \cross \E$ as $M$ with $M_1(e,e') \le M(e,e')$
and with strict inequality for at least one pair $(e,e') \in\E \cross \E$, and since $\lambda$ is the
Perron-Frobenius eigenvalue of $M$, this inequality follows from Perron-Frobenius theory.  See for
example  Theorem~4.4.7 and Theorem~4.4.4 of \cite{LindMarcus:symbolic}. 
\end{proof}

Now we relate Fact~\ref{FactLambdaPrime} to the asymptotic compression factor of $k$ acting on
leaves of $\Hull\fol^s$. 

\begin{fact}
\label{FactAsymptoticCompression}
For each $\epsilon>0$ there exists an integer $N \ge 1$ such that for all $n \ge N$, if $\alpha$ is an arc
in a leaf $\ell$ of $\Hull\fol^s$, and if $\Length_K(\alpha)$ is sufficiently long (depending on $n$),
then 
$$\frac{\Length_K(\alpha)}{\Length_K(k^n \alpha)} \le (\lambda' + \epsilon)^n 
$$
\end{fact}

\begin{proof} Applying Fact~\ref{FactHullInfo}, for each leaf $\ell$ of $\Hull\fol^s$ let $\bar
k(\ell)$ denote the unique leaf of $\Hull\fol^s$ such that $k(\bar k(\ell))=\ell$, so $\bar k(\ell)
\subset k^\inv(\ell)$. We extend the set map $\bar k$ to subsets of leaves of $\Hull\fol^s$, by setting
$\bar k(\alpha) = \bar k(\ell) \intersect k^\inv(\alpha)$ for each $\alpha \subset\ell$.

Let $\ell$ be a leaf of $\Hull\fol^s$ and $\alpha\subset\ell$ an arc. Let $k^n \alpha \intersect G =
\{x_0,x_1,\ldots,x_I\}$ in order. Note that if $0<i<I$ then $\bar k^n(x_i)
\subset \alpha$, so we can augment $\alpha$, replacing it by $\bar k^n(x_0) \union \alpha \union \bar
k^n(x_I)$, and increasing the length of $\alpha$ by an amount depending only on $n$, without changing $k^n
\alpha$. If $\Length_K(\alpha)$ were sufficiently long to start with then this increase would change
$\Length_K(\alpha)$ by an arbitrarily small factor. We may therefore assume that $\alpha = \bar k^n(k^n
\alpha)$. 

It follows that $\alpha$ is equal to the disjoint union of the arcs $\bar k^n(x_0),\ldots,\bar k^n(x_I)$
together with $I$ additional open vertical wedge segments, each of which maps homeomorphically by $k^n$
to the $I$ open wedge segments that constitute the arc $k^n \alpha$. We therefore have
$$\Length_K(\alpha) = \sum_{i=0}^I \Length_K(\bar k^n(x_i)) + I
$$
and dividing by $\Length_K(k^n \alpha)=I$ we get
$$\frac{\Length_K(\alpha)}{\Length_K(k^n \alpha)} = \frac{1}{I} \sum_{i=0}^I \Length_K(\bar k^n(x_i)) + 1
$$
We also have $\Length_K(\bar k^n(x_i)) \le \abs{g^{-n}(x_i) \intersect I_1}$. By combining this with
Fact~\ref{FactLambdaPrime} it follows that if $n$ is sufficiently large then $\Length_K(\bar k^n(x_i))
\le (\lambda'+\frac{\epsilon}{2})^n$ and so
\begin{align*}
\frac{\Length_K(\alpha)}{\Length_K(k^n \alpha)} &\le \left(\lambda'+\frac{\epsilon}{2}\right)^n +
\frac{\lambda'+\frac{\epsilon}{2}}{I} + 1 \\
\intertext{By taking $n$ sufficiently large, this last quantity is}
  &\le (\lambda'+\epsilon)^n
\end{align*}

\end{proof}

Now we prove Theorem~\ref{TheoremLambdaBigger}. Choose a train track representative $\gamma \from \Gamma \to \Gamma$ of $\phi^\inv$ and let $\tilde \gamma \from \wt \Gamma \to \wt \Gamma$ be a lift to the universal cover. Since $g \from G \to G$ is a train track representative of $\phi$, we can choose the lifts $\tilde\gamma \from\tilde\Gamma \to \tilde\Gamma$ and $\tilde g \from \tilde G \to \tilde G$ to represent inverse automorphisms of $F_r$. It follows that there is an $F_r$-equivariant quasi-isometry $h\from \wt G\to \wt\Gamma$ such that $h$ is a ``quasiconjugacy'' between $\tilde k$ and ``$\tilde \gamma^\inv$'', meaning that $d(\tilde \gamma \composed h \composed \tilde g(x),h(x))$ is uniformly bounded over all $x\in\wt G$. Now we apply this to a particular leaf of $\wt\Lambda^u(\phi^\inv)$ as follows. Pick a leaf $\ell'$ of $\Lambda^u(\phi^\inv)$ realized in~$\Gamma$. For convenience we assume $\ell'$ is periodic under $\gamma$. By passing to a power of $\phi$ we may assume that $\ell'$ is fixed by $\gamma$. We may choose
a lift $\tilde\ell' \subset \wt\Gamma$, and we may choose the lifts $\tilde g$ and $\tilde \gamma$, so that $\tilde\ell'$ is fixed by $\tilde\gamma$. Let $\tilde\ell$ denote the corresponding leaf of $\Hull\wt\fol^s$, and so $\tilde\ell$ is fixed by $\tilde k$. Pushing $\tilde\ell$ into $\wt G$, and mapping over by $h$ to $\tilde \Gamma$, the result is Hausdorff equivalent to $\tilde\ell'$, and composing by the closest point projection to $\tilde\ell'$, we obtain a map still denoted $h \from \tilde\ell\to\tilde\ell'$ which is a quasi-isometry and a quasiconjugacy, that is, there are constants $\kappa\ge 1$, $\eta \ge 0$ such that
$$\frac{1}{\kappa} \, d_\ell(x,y) - \eta \le d_{\ell'}(hx,hy) \le \kappa \, d_\ell(x,y) + \eta,
\quad\text{for all $x,y \in \tilde\ell$}
$$
and such that
$$d_{\ell'}\left(\tilde \gamma \composed h \composed \tilde k(x),h(x)\right) \le \eta \quad\text{for all
$x\in\tilde\ell$.}
$$
Here we use $d_\ell$ for distance in $\tilde\ell$ and $d_{\ell'}$ for distance in $\tilde\ell'$.
By induction we obtain for each $n$ a constant $\eta_n \ge 0$ such that
$$d_{\ell'}\left(\tilde\gamma^n \composed  h \composed  \tilde k^n(x), h(x)\right) \le \eta_n,
\quad\text{for all $x \in \ell'$}
$$
We know that $\tilde \gamma$ expands length on $\tilde\ell'$ by the exact factor of $\lambda(\phi^\inv)$. On
the other hand, Fact~\ref{FactAsymptoticCompression} tells us that $\tilde k$ contracts length on
$\tilde\ell$ by an asymptotic factor of at most~$\lambda'$. Combining these, we now show that
$\lambda(\phi^\inv) \le \lambda'$. 


By Fact~\ref{FactAsymptoticCompression}, for each sufficiently long leaf segment $[x,y]$ of
$\tilde\ell$ we have
\begin{align*}
d_\ell(x,y) &\le (\lambda'+\epsilon)^n \, d_\ell(\tilde k^n (x), \tilde k^n (y)) \\
\intertext{and we also have}
d_\ell(\tilde k^n (x), \tilde k^n (y)) &\le \kappa \, d_{\ell'}(h \composed \tilde k^n (x), h \composed \tilde
k^n (y)) + \kappa\eta \\
d_{\ell'}(h \composed \tilde k^n (x), h \composed \tilde k^n (y)) &= \lambda(\phi^\inv)^{-n}
d_{\ell'}(\tilde\gamma^n \composed h \composed \tilde k^n (x), \tilde \gamma^n \composed h \composed \tilde
k^n (y)) \\
d_{\ell'}(\tilde\gamma^n \composed h \composed \tilde k^n (x), \tilde \gamma^n \composed h \composed \tilde
k^n (y)) &\le d_{\ell'}(hx,hy) + 2 \eta_n \\
d_{\ell'}(hx,hy) &\le \kappa \, d_\ell(x,y) + \eta 
\end{align*}
Combining these we get
\begin{align*}
d_\ell(x,y) &\le (\lambda'+\epsilon)^n \lambda(\phi^\inv)^{-n} \kappa^2 d_\ell(x,y) +
\underbrace{(\lambda'+\epsilon)^n(\eta + 2 \eta_n + \kappa\eta)}_{B_n} \\
1 &\le (\lambda'+\epsilon)^n \lambda(\phi^\inv)^{-n} \kappa^2 + \frac{B_n}{d_\ell(x,y)}
\end{align*}
By Fact~\ref{FactAsymptoticCompression}, we can let $d_\ell(x,y) \to \infinity$ and we get
\begin{align*}
\lambda(\phi^\inv)^n &\le (\lambda'+\epsilon)^n \kappa^2 \\
\lambda(\phi^\inv) &\le (\lambda'+\epsilon) \kappa^{2/n} \\
\intertext{and then we can let $n \to \infinity$ to get}
\lambda(\phi^\inv) &\le \lambda'+\epsilon \\
\intertext{Finally, letting $\epsilon \to 0$ we get}
\lambda(\phi^\inv) &\le \lambda' < \lambda(\phi)
\end{align*}
This completes the proof of Theorem~\ref{TheoremLambdaBigger}.


\newcommand{\etalchar}[1]{$^{#1}$}
\providecommand{\bysame}{\leavevmode\hbox to3em{\hrulefill}\thinspace}
\providecommand{\MR}{\relax\ifhmode\unskip\space\fi MR }
\providecommand{\MRhref}[2]{%
  \href{http://www.ams.org/mathscinet-getitem?mr=#1}{#2}
}
\providecommand{\href}[2]{#2}

\bigskip

\bigskip\noindent
\textsc{\scriptsize
Michael Handel:\\
Department of Mathematics and Computer Science\\
Lehman College - CUNY\\
250 Bedford Park Boulevard W\\
Bronx, NY 10468\\
michael.handel@lehman.cuny.edu
}

\bigskip

\bigskip\noindent
\textsc{\scriptsize
Lee Mosher:\\
Department of Mathematics and Computer Science\\
Rutgers University at Newark\\
Newark, NJ 07102\\
mosher@andromeda.rutgers.edu
}

\end{document}